\documentclass[a4paper,12pt]{article}
\usepackage{amsmath}
\usepackage{amsthm}
\usepackage{amssymb}
\usepackage{amscd}
\usepackage{enumerate}
\theoremstyle{definition}
\numberwithin{equation}{section}
\newtheorem{thm}{Theorem}[section]
\newtheorem{dfn}[thm]{Definition}
\newtheorem{exa}[thm]{Example}
\newtheorem{prop}[thm]{Proposition}
\newtheorem{cor}[thm]{Corollary}
\newtheorem{lem}[thm]{Lemma}
\newtheorem{conj}[thm]{Conjecture}
\newtheorem{rem}[thm]{Remark}

\def\Ker{\mathop{\mathrm{Ker}}\nolimits}

\def\Ker{\mathop{\mathrm{Ker}}\nolimits}

\def\mod{\mathop{\mathrm{mod}}\nolimits}

\def\End{\mathop{\mathrm{End}}\nolimits}

\def\Ind{\mathop{\mathrm{Ind}}\nolimits}

\def\Supp{\mathop{\mathrm{Supp}}\nolimits}
\def\row{\mathop{\mathrm{row}}\nolimits}
\def\col{\mathop{\mathrm{col}}\nolimits}
\def\res{\mathop{\mathrm{res}}\nolimits}
\def\vac{\mathop{\mathrm{vac}}\nolimits}
\def\triv{\mathop{\mathrm{triv}}\nolimits}

\def\rad{\mathop{\mathrm{rad}}\nolimits}

\newcommand{\mf}[1]{{\mathfrak{#1}}}

\newcommand{\bb}[1]{{\mathbb{#1}}}
\title{\textbf{Composition Factors of Polynomial Representation of DAHA and
Crystallized Decomposition Numbers}}
\author{Naoya ENOMOTO \\
henon@kurims.kyoto-u.ac.jp}
\date{Resarch Institute of Mathematical Sciences, Kyoto University}
\begin{document}
\maketitle
\begin{abstract}
We determine the composition factors of the polynomial representation of DAHA,
conjectured by M. Kasatani in \cite[Conjecture 6.4.]{Kasa}. He constructed
an increasing sequence of subrepresentations in the polynomial
representation of DAHA using the ``multi-wheel condition'', and conjectured that it is a composition series. On the other hand, 
DAHA has two degenerate versions called the ``degenerate DAHA'' and the ``rational DAHA''. The category $\mathcal{O}$ of modules over these three algebras and the category of modules over the $v$-Schur algebra are closely related. By using this relationship, We reduce
the determination of composition factors of polynomial representations of 
DAHA to the determination of the composition factors of the Weyl module
$W^{(1^n)}$ for the $v$-Schur algebra. By using the LLT-Ariki type theorem of
$v$-Schur algebra proved by Varagnolo-Vasserot, we determine the
composition
factors of $W^{(1^n)}$ by calculating the upper global basis and crystal
basis of Fock space of $U_q(\widehat{\mf{sl}}_\ell)$.
\end{abstract}

\textbf{Acknowledgement.} I would like to thank Professor Masaki Kashiwara for
his many advises and suggestions, especially about the global basis, and
Professor Takeshi Suzuki for his advises and suggestions about the double affine
Hecke algebra and its degenerate versions, especially the reduction of the
determination of composition facters. I also would like to thank Masahiro
Kasatani, Toshiro Kuwabara for their suggestions and useful discussions about representations of double
affine Hecke algebra and its degenerate versions, and Professor Susumu Ariki for his
suggestions about decomposition numbers and their combinatorics. I also would like to thank Professor Hyoe Miyachi for his suggestions about his prior works and related topics.
\newpage

\tableofcontents
\newpage

\section{Introduction}
\paragraph{Double Affine Hecke Algebra}
The double affine Hecke algebra (DAHA) is a 2-parameter analogue of the Iwahori-Hecke algebra introduced by I. Cherednik in \cite{Ch1}. This algebra is closely related to the symmetric or orthogonal polynomials. In \cite{Ch2}, Cherednik proved the Macdonald inner product conjecture by using DAHA. \\
\indent The DAHA $\mathcal{H}_n$ of type $GL_n$ is generated by 
\[
T_i \ (1 \le i \le n-1), X_j^{\pm{1}}, Y_j^{\pm{1}} \ (1 \le j \le n),
\]
and has two parameters $\zeta$ and $\tau$. The generators $T_i$ satisfy the Hecke relation $(T_i-\zeta^{1/2})(T_i+\zeta^{-1/2})=0$ and generate the Iwahori-Hecke algebra of type $A$. The two subalgebras $\langle T_i (1 \le i \le n-1), X_j (1 \le j \le n) \rangle$ and $\langle T_i (1 \le i \le n-1), Y_j (1 \le j \le n)$ are both isomorphic to the affine Hecke algebra of type $GL_n$. The parameter $\tau$ appears in some relations between $X$ and $Y$. \\
\indent The DAHA has a faithful representation on $\bb{C}(\zeta^{1/2},\tau)[X_1^{\pm{1}}, \cdots X_n^{\pm{1}}]$ defined by difference Dunkl operators. This representation is called the polynomial representation. If $\zeta$ and $\tau$ are generic, the non-symmetric Macdonald polynomials are simultaneous $Y$-eigenvectors. \\
\indent The category consisting of $Y$-locally finite modules of DAHA is called the ``category $\mathcal{O}$''. Especially, the combinatorial construction of $Y$-semisimple irreducible representations is obtained by T. Suzuki and M. Vazirani in \cite{SV}. On the other hand, in \cite{Va}, E. Vasserot studied the geometric construction of the DAHA and its irreducible representations and some composition multiplicities. \\
\paragraph{The Kasatani Conjecture}
Recall that the DAHA has a faithful representation on the Laurent polynomial ring, called the polynomial representation. If $\zeta$ and $\tau$ are generic, this representation is irreducible and $Y$-semisimple. But if the two parameters of DAHA specialized at $\zeta^{\ell}\tau^{r}=1$, then the polynomial representation is not any more irreducible and $Y$-semisimple in general. \\
\indent M. Kasatani constructed in \cite{Kasa} the increasing sequence of the subrepresentation of the polynomial representation by using the ``multi-Wheel condition''. He conjectured that this sequence is a composition series of the polynomial representation of DAHA. 
\paragraph{Main results}
In this paper we prove Kasatani's conjecture when
\[
(\ell,r)=1 \ \text{and} \ \ell \neq 2. 
\]
\quad {} \\
\quad {} \\
\indent In the rest of this introduction, we explain the strategy of our proof of this main results. \\
\quad {} \\
\indent The DAHA has two degenerate versions, simply called the ``degenerate DAHA'' and the ``rational DAHA''. 
\paragraph{Degenerate DAHA}
Roughly speaking the degenerate DAHA $\bb{H}^{deg}_{n,h}$ is obtained from DAHA by degenerating $T$ and $Y$. This algebra is generated by the following three subalgebras:
\[
\bb{C}[X_1^{\pm{1}}, \cdots ,X_n^{\pm{1}}], \ \bb{C}\mf{S}_n, \ \bb{C}[y_1, \cdots ,y_n],
\]
and this algebra has one parameter $h$. \\
\indent The degenerate DAHA has the category $\bb{O}^{deg}_h$ consisting of $y$-locally finite modules. By T. Suzuki in \cite{Su}, \cite{Su2}, this category and the irreducible modules are studied. Especially this category has the standard modules induced from the irreducible representations of the symmetric group. \\
\indent In the aspect of the relation with the conformal field theory, see in \cite{AST}.
\paragraph{Rational DAHA}
Roughly speaking the rational DAHA $\bb{H}^{rat}_{n,h}$ is obtained from DAHA by degenerating $T, X$ and $Y$. This algebra is generated by the following three subalgebras:
\[
\bb{C}[x_1, \cdots ,x_n], \ \bb{C}\mf{S}_n, \ \bb{C}[y_1, \cdots ,y_n],
\]
and this algebra has one parameter $h$. \\
\indent The rational DAHA has also the category $\bb{O}^{rat}$ consisting of $y$-locally nilpotent modules. Especially this category has the standard modules induced from the irreducible representations of the symmetric group. This algebra is a special version of the symplectic reflection algebra defined by P. Etingof and V. Ginzburg in \cite{EG}. The category $\bb{O}^{rat}_h$ are studied in \cite{GGOR}. \\
\indent In the aspect of the relation with the coherent sheaves on the Hilbert schemes, see \cite{GS1}, \cite{GS2}.
\paragraph{$v$-Schur algebra}
The $v$-Schur algebra $\bb{S}(n)$ has two different definitions. One is the endomorphism ring of the permutation module of the Iwahori-Hecke algebra, introduced by R. Dipper and G. James in \cite{DJ}. Second is the quotient of the quantum universal enveloping algebra $U_v(\mf{gl}_n)$ in the tensor representation of the vector representation, introduced by A. Beilinson, G. Lusztig, R. MacPherson in \cite{BLM}. \\
\indent This algebra is a cellular algebra introduced by J. Graham and G. Lehrer in \cite{GL}. Moreover this is a quasi-hereditary cover of the Iwahori-Hecke algebra, see in \cite{Rou2}. Thus the ordinary and modular representation theory of the $v$-Schur algebra are studied well. For example, the complete collection of the ordinary irreducible representations are indexed by partitions of $n$, called the Weyl modules $W^{\lambda}$. If $v$ is a root of unity, the complete collection of irreducible representations are indexed by partitions of $n$, denoted by $L^\lambda:=W^\lambda/\rad{W^\lambda}$.\\
\paragraph{Relationships of these algebras}
The category $\mathcal{O}$ of DAHA and its two degenerate versions and the category of the $v$-Schur algebra are related as seen in the following figure;
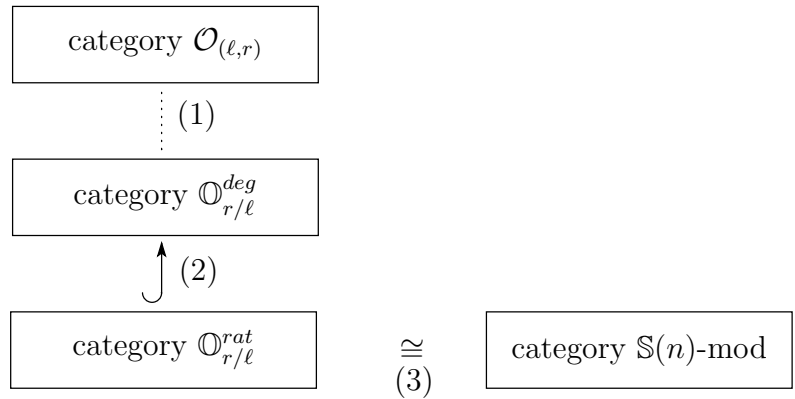
\[
\unitlength 0.1in
\begin{picture}(49.90,20.00)(0.90,-22.00)
\put(18.0000,-4.0000){\makebox(0,0){category
$\mathcal{O}_{(\ell,r)}$}}%
%
\special{pn 8}%
\special{pa 1000 200}%
\special{pa 2600 200}%
\special{pa 2600 600}%
\special{pa 1000 600}%
\special{pa 1000 200}%
\special{fp}%
\put(18.0000,-12.0000){\makebox(0,0){category $\bb{O}_{r/\ell}^{deg}$}}%
%
\special{pn 8}%
\special{pa 1000 1000}%
\special{pa 2600 1000}%
\special{pa 2600 1400}%
\special{pa 1000 1400}%
\special{pa 1000 1000}%
\special{fp}%
\put(42.8000,-20.0000){\makebox(0,0){category $\bb{S}(n)$-mod}}%
%
\special{pn 8}%
\special{pa 3480 1800}%
\special{pa 5080 1800}%
\special{pa 5080 2200}%
\special{pa 3480 2200}%
\special{pa 3480 1800}%
\special{fp}%
\put(17.9000,-20.0000){\makebox(0,0){category $\bb{O}_{r/\ell}^{rat}$}}%
%
\special{pn 8}%
\special{pa 990 1800}%
\special{pa 2590 1800}%
\special{pa 2590 2200}%
\special{pa 990 2200}%
\special{pa 990 1800}%
\special{fp}%
%
\special{pn 8}%
\special{pa 1780 650}%
\special{pa 1780 950}%
\special{dt 0.045}%
\special{pa 1780 950}%
\special{pa 1780 949}%
\special{dt 0.045}%
\put(30.9000,-20.0000){\makebox(0,0){$\cong$}}%
%
\special{pn 8}%
\special{pa 1780 1700}%
\special{pa 1780 1440}%
\special{fp}%
\special{sh 1}%
\special{pa 1780 1440}%
\special{pa 1760 1507}%
\special{pa 1780 1493}%
\special{pa 1800 1507}%
\special{pa 1780 1440}%
\special{fp}%
%
\special{pn 8}%
\special{ar 1730 1700 50 50  6.2831853 6.2831853}%
\special{ar 1730 1700 50 50  0.0000000 3.1415927}%
\put(18.6000,-6.9000){\makebox(0,0)[lt]{(1)}}%
\put(18.7000,-15.1000){\makebox(0,0)[lt]{(2)}}%
\put(31.0000,-21.7000){\makebox(0,0){(3)}}%
\end{picture}%
\]
\begin{quote}
(1) \ \textbf{Varagnolo-Vasserot's block equivalence \cite{VV2}} \\
The category $\mathcal{O}$ of the DAHA and the category of $\bb{O}^{deg}$ are not equivalent. But we can consider the specialized DAHA $\mathcal{H}_n^{(\ell.r)}$ where the two parameters specialized at $\zeta^{\ell}\tau^{r}=1$ and $(\zeta,\tau)=1$. Let $\mathcal{O}_{(\ell,r)}$ be the category $\mathcal{O}$ of $\mathcal{H}_n^{(\ell,r)}$. Further, we can consider the subcategory ${}^{\chi}\mathcal{O}_{(\ell,r)}$ of $\mathcal{O}_{(\ell,r)}$ consisting the modules such that all the $Y$-weights of them belong to the affine Weyl group orbit of $\chi$. And we can consider the similar full subcategory ${}^{\chi}\bb{O}^{deg}$ of $\bb{O}^{deg}$. If $\chi$ satisfies some conditions, then the categories ${}^{\chi}\mathcal{O}_{(\ell,r)}$ and ${}^{\chi}\bb{O}_{r/\ell}$ are equivalent. This equivalence is a direct generalization of the equivalence between some category of representations of the affine Hecke algebra and the one of the degenerate affine Hecke algebra proved by G. Lusztig in \cite{Lus}. \\
\quad {} \\
(2) \ \textbf{T. Suzuki's embedding \cite{Su}} \\
T. Suzuki proved that the rational DAHA can be embedded in the degenerate DAHA. Moreover the functer $\bb{H}^{deg} \otimes_{\bb{H}^{rat}}-$ is fully faithful and exact. Thus the standard module and its (unique) simple quotient are sent to the standard module and its (unique) simple quotient under the above functer.\\
\quad {} \\
(3) \ \textbf{R. Rouquier's equivalence \cite{Rou2}} \\
R. Rouquier proved that the category $\bb{O}^{rat}_{r/\ell}$ is equivalent to the category of $\bb{S}(n)$-mod at $v=\sqrt[\ell]{1}$ unless $\ell \neq 2$. And he proved that the standard modules are sent to the Weyl modules of the $v$-Schur algebra. This had been conjectured in \cite{GGOR}. \\
\quad {} \\
(4) \ \textbf{Lascoux-Leclerc-Thibon's conjecture} \\
Since S. Ariki proved the LLT conjecture on the decomposition numbers of the cyclotomic Hecke algebra in \cite{Ari1}, the modular representation theory of Hecke algebras are closely related to the representation theory of the quantum enveloping algebra and the theory of the crystal and global basis. Moreover M. Varagnolo and E. Vasserot proved the extended version of the LLT conjecture about the decomposition numbers of $v$-Schur algebra in \cite{VV1}. By this result, the decomposition numbers $[W^{\lambda}:L^{\mu}]$ for $v$-Schur algebras are described by the transition matrix of the standard basis and the global basis of the Fock space of $U_q(\mf{sl}_\ell)$.
\end{quote}
Therefore we can reduce some problems of the representation theory of DAHA to some calculation of the global basis of the Fock space. The main result of this paper is based on this strategy.

By using the above strategy, we can reduce the determination of the composition factors of the polynomial representations to the one of the Weyl module $W^{(1^n)}$. It is again equivalent to the determination of the coefficient of the upper global basis $G^{up}(\mu)$ in the expansion by standard basis $|(n)\rangle$ in the Fock space. We calculate this coefficient and obtain the following
\[
[W^{(1^n)}:L^{\mu}]=\left\{
\begin{array}{ll}
1 & \text{if} \ \mu'=\mu_i^{(n)} (0 \le i \le N=[r/\ell])\\
0 & \text{otherwise} 
\end{array}
\right..
\]
Here $\mu'$ is the conjugate partition of $\mu$. On the definition of the partition $\mu_i^{(n)}$, see Definition \ref{defp}. Note that this result is true for $\ell=2$. However since Rouquier's equivalence (3) is not proved at $\ell=2$, we cannot prove the Kasatani conjecture at $\ell=2$ by our method. 
\newpage

\section{DAHA and the Kasatani Conjecture}
\setcounter{subsection}{-1}
\subsection{Notations for affine root systems and affine Weyl groups}
We will use the following notations for the affine root system and the affine Weyl group of type $A$. \\
\indent Let $\mf{h}$ be an $(n+2)$-dimensional vector space over $\bb{C}$ with basis
\[
\mf{h}=\bigoplus_{i=1}^{n}\bb{C}\varepsilon_i^{\vee} \oplus \bb{C}c \oplus \bb{C}d.
\]
Let $\mf{h}^*$ be the dual space of $\mf{h}$, where $\varepsilon_i,\Lambda$ and $\delta$ are the dual basis of $\varepsilon_i^{\vee},c$ and $d$;
\[
\mf{h}^*=\bigoplus_{i=1}^{n}\bb{C}\varepsilon_i \oplus \bb{C}\Lambda \oplus \bb{C}\delta.
\]
There exists a non-degenerate symmetric bilinear form $( \ | \ )$ on $\mf{h}$ defined by
\[
(\varepsilon_i^{\vee}|\varepsilon_j^{\vee})=\delta_{ij}, \quad
(\varepsilon_i^{\vee}|c)=(\varepsilon_i^{\vee}|d)=0, \quad
(c|d)=1, \quad 
(c|c)=(d|d)=0.
\]
The natural pairing is denoted by $\langle \ | \ \rangle:\mf{h}^* \times \mf{h} \to \bb{C}$. There exists an isomorphism $\mf{h}^* \to \mf{h}$ such that 
\[
\varepsilon_i \mapsto \varepsilon_i^{\vee}, \quad \delta \mapsto c, \quad \Lambda \mapsto d.
\]
We denote by $h^{\vee} \in \mf{h}$ the image of $h \in \mf{h}^*$ under this isomorphism. We can introduce the bilinear form $( \ | \ )$ on $\mf{h}^*$ through this isomorphism, and then
\[
(h|k)=\langle{h}|k^{\vee}\rangle=(h^{\vee}|k^{\vee}) \ \text{for} \ h,k \in \mf{h}^*.
\]
\indent Put 
\[
\alpha_{ij}=\varepsilon_i-\varepsilon_j, \ (1 \le i \neq j \le n) \quad \alpha_i=\alpha_{i,i+1} \ (1 \le i \le n-1).
\]
Then 
\begin{eqnarray*}
R&=&\{\alpha_{ij}|1 \le i \neq j \le n\} \subset \mf{h}^*, \\
R^+&=&\{\alpha_{ij} \in R|i<j\}, \\
\Pi&=&\{\alpha_1, \cdots ,\alpha_{n-1}\}
\end{eqnarray*}
give the set of roots, positive roots and simple roots of type $A_{n-1}$, respectively. \\
\indent Put 
\[
\alpha_0=-\alpha_{1n}+\delta.
\]
Then
\begin{eqnarray*}
\widehat{R}&=&\{\alpha+k\delta|\alpha \in R, \ k \in \bb{Z}\} \subset \mf{h}^*, \\
\widehat{R}^+&=&\{\alpha+k\delta|\alpha \in R^+, k \in \bb{Z}_{\ge 0}\} \sqcup \{-\alpha+k\delta|\alpha \in R^+, \ k \in \bb{Z}_{\ge 1}\}, \\
\widehat{\Pi}&=&\{\alpha_0,\alpha_1, \cdots ,\alpha_{n-1}\}
\end{eqnarray*}
give the set of real roots, positive real roots and simple roots of type $A_{n-1}^{(1)}$, respectively. \\
\indent Let $P$ and $P^{\vee}$ be the weight lattice and co-weight lattice defined by
\[
P=\bigoplus_{i=1}^{n}\bb{Z}\varepsilon_i \subset \mf{h}^{*}, \quad P^{\vee}=\bigoplus_{i=1}^{n}\bb{Z}\varepsilon_i^{\vee} \subset \mf{h}.
\]
\indent We introduce the affine Weyl group of type $A_{n-1}^{(1)}$.
\begin{dfn}
\quad The extended affine Weyl group $W_n$ of type $A_{n-1}^{(1)}$ is the group defined by the following generators and relations;
\begin{eqnarray*}
\text{generators} &:& s_0,s_1, \cdots ,s_{n-1},\pi^{\pm{1}}, \\
\text{relations} &:& s_i^2=1 \quad (0 \le i \le n-1), \\
&{}&s_is_{i+1}s_i=s_{i+1}s_is_{i+1} \quad (i \in \bb{Z}/n\bb{Z}, n>2), \\
&{}&s_is_j=s_js_i \quad (j \not\equiv i,i\pm{1}), \\
&{}&{\pi}s_i=s_{i+1}\pi \quad (i \in \bb{Z}/n\bb{Z}), \\
&{}&\pi^{-1}\pi=\pi\pi^{-1}=1. \\
\end{eqnarray*}
The subgroup $\langle s_1, \cdots ,s_{n-1}\rangle$ is isomorphic to the symmetric groups $\mf{S}_n$. The subgroup $\langle s_0,s_1, \cdots ,s_{n-1}\rangle$ is called the (non-extended) affine Weyl group of type $A_{n-1}^{(1)}$.
\end{dfn}
We can describe the extended affine Weyl group as a semi-direct product group. Put
\[
X_{\varepsilon_1}={\pi}s_{n-1}s_{n-2} \cdots s_1, \quad X_{\varepsilon_i}=\pi^{i-1}X_{\varepsilon_1}\pi^{-i+1} \ (2 \le i \le n).
\]
Then there exists an embedding $P \hookrightarrow W_n$ defined by $\varepsilon_i \mapsto X_{\varepsilon_i}$. We denote $X_{\eta}$ the image of $\eta \in P$ under this embedding. Then there exists an isomorphism
\[
W_n \cong P \rtimes \mf{S}_n=\langle{X}_{\varepsilon_i} \ (1 \le i \le n),s_1, \cdots ,s_{n-1}\rangle
\]
such that 
\[
wX_{\eta}w^{-1}=X_{w(\eta)} \ (w \in \mf{S}_n, \ \eta \in P).
\]
Here the symmetric group $\mf{S}_n$ acts on $P$ by $s_i:\varepsilon_i \leftrightarrow \varepsilon_{i+1}$.
\indent The extended affine Weyl group $W_n$ acts on $\mf{h}^*$ by
\begin{eqnarray*}
s_i(h)&=&h-(\alpha_i|h)\alpha_i \quad (h \in \mf{h}^*), \\
X_{\eta}(h)&=&h+(\delta|h)\eta-\left\{(\eta|h)+\frac{1}{2}(\eta|\eta)(\delta|h)\right\}\delta, \\
\pi(\varepsilon_i)&=&\varepsilon_{i+1} \quad (1 \le i \le n-1), \\
\pi(\varepsilon_n)&=&\varepsilon_n-\delta, \\
\pi(\Lambda)&=&\Lambda, \\
\pi(\delta)&=&\delta.
\end{eqnarray*}
The dual action on $\mf{h}$ is the following;
\begin{eqnarray*}
s_i(h^{\vee})&=&h^{\vee}-\langle\alpha_i|h^{\vee}\rangle\alpha_i^{\vee} \quad (h^{\vee} \in \mf{h}), \\
X_{\eta}(h^{\vee})&=&h^{\vee}+\langle\delta|h^{\vee}\rangle\eta^{\vee}-\left\{\langle\eta|h^{\vee}\rangle+\frac{1}{2}(\eta|\eta)\langle\delta|h^{\vee}\rangle\right\}c, \\
\pi(\varepsilon_i^{\vee})&=&\varepsilon_{i+1}^{\vee} \quad (1 \le i \le n-1), \\
\pi(\varepsilon_n^{\vee})&=&\varepsilon_n^{\vee}-c, \\
\pi(c)&=&c, \\
\pi(d)&=&d.
\end{eqnarray*}
\subsection{Double affine Hecke algebra and its Polynomial representation}
\subsubsection{Double affine Hecke algebra of type $GL_n$}
Let $\bb{K}$ be a field $\bb{C}(\zeta^{1/2},\tau)$. The double affine Hecke
algebra of type $GL_n$ is defined as follows.
\begin{dfn} The double affine Hecke algebra $\mathcal{H}_n$ of type $GL_n$
is an associative algebra over $\bb{K}$ generated by
\[
T_i \ (0 \le i \le n-1), \quad Y_{\eta} \ (\eta \in P \oplus \bb{Z}\delta), \quad \pi^{\pm{1}}
\]
satisfying the following relations
\begin{eqnarray*}
\begin{array}{ll}
Y_{\delta}=\tau, & \\
(T_i-\zeta^{1/2})(T_i+\zeta^{-1/2})=0 & (0 \le i \le n-1), \\
T_iT_{i+1}T_i=T_{i+1}T_iT_{i+1} & (i \in \bb{Z}/n\bb{Z}), \\
T_iT_j=T_jT_i & (\text{otherwise}), \\
T_iY_{\eta}-Y_{s_i(\eta)}T_i=(\zeta^{1/2}-\zeta^{-1/2})\dfrac{Y_{s_i\eta}-Y_{\eta}}{Y_{\alpha_i}-1} & (0 \le i \le n-1), \\
\pi{T_i}=T_{i+1}\pi & (i \in \bb{Z}/n\bb{Z}), \\
\pi{Y}_{\eta}=Y_{\pi(\eta)}\pi, & \\
Y_{\eta}Y_{\xi}=Y_{\eta+\xi}. & 
\end{array}
\end{eqnarray*}
The two subalgebras
\begin{eqnarray*}
H_n&=&\langle T_1, \cdots ,T_{n-1} \rangle, \\
H_n^{aff}&=&\langle T_1, \cdots ,T_{n-1},Y_\eta \ (\eta \in P)\rangle
\end{eqnarray*}
are isomorphic to the Iwahori-Hecke algebra of type $A_{n-1}$ and the affine Hecke algebra of type $GL_n$, respectively.
\end{dfn}
\begin{rem}
Put $Y_i=Y_{\varepsilon_i}$ and 
\[
X_1=T_1 \cdots T_{n-1}\pi^{-1}, \quad X_i=\pi^{i-1}X_1\pi^{-i+1}.
\]
There is an another description of generators and relations as the following;
\begin{eqnarray*}
\begin{array}{lll}
\text{generators}: & T_i \ (1 \le i \le n-1), \quad Y_j^{\pm{1}}, X_j^{\pm{1}} \ (1 \le j \le n), & \\
\text{relations}: & (T_i-\zeta^{1/2})(T_i+\zeta^{-1/2})=0 & (1 \le i \le n-1), \\
& T_iT_{i+1}T_i=T_{i+1}T_iT_{i+1} & (1 \le i \le n-1), \\
& T_iT_j=T_jT_i & (|i-j| \ge 2), \\
& T_iX_{i+1}T_i=X_i & (1 \le i \le n-1), \\
& T_iX_j=X_jT_i & (j \neq i,i+1), \\
& T_iY_iT_i=Y_{i+1} & (1 \le i \le n-1), \\
& T_iY_j=Y_jT_i & (j \neq i,i+1), \\
& X_2^{-1}Y_1X_2Y_1^{-1}=T_1^2, & \\
& \displaystyle X_j\left(\prod_{k=1}^{n}Y_k\right)=\tau
\left(\prod_{k=1}^{n}Y_k\right)X_j & (1 \le j \le n), \\
& \displaystyle Y_j\left(\prod_{k=1}^{n}X_k\right)=\tau
\left(\prod_{k=1}^{n}X_k\right)Y_j & (1 \le j \le n), \\
& X_iX_j=X_jX_i, X_iX_i^{-1}=1 & (1 \le i,j \le n), \\
& Y_iY_j=Y_jY_i, Y_iY_i^{-1}=1 & (1 \le i,j \le n).
\end{array}
\end{eqnarray*}
\end{rem}

\subsubsection{Polynomial representation}
The DAHA $\mathcal{H}_n$ has a representation on $\bb{K}[x_1^{\pm{1}}, \cdots ,x_n^{\pm{1}}]$ defined by the difference Dunkl operators. This representation is called the polynomial representation. 
\begin{prop} \quad {} \\
(1) \ The DAHA $\mathcal{H}_n$ has a representation on $\bb{K}[x_1^{\pm{1}}, \cdots ,x_n^{\pm{1}}]$ defined by the following,
\begin{eqnarray*}
X_j &\mapsto& x_j \ (\text{multiplication}), \\
T_i &\mapsto& \zeta^{1/2}s_i+\frac{\zeta^{1/2}-\zeta^{-1/2}}{x_{i+1}x_i^{-1}-1}(s_i-1), \\
Y_j &\mapsto& T_j^{-1} \cdots T_{n-1}^{-1}{\omega}T_1 \cdots T_{j-1},
\end{eqnarray*}
where $s_i$ is the permutation of $x_i$ and $x_{i+1}$, and 
\[
(\omega{f})(x_1, \cdots ,x_n)=f(\tau^{-1}x_n, x_1, \cdots ,x_{n-1}).
\]
(Note that $\pi=X_1^{-1}T_1 \cdots T_{n-1}$.)\\
(2) \ This representation is the induced representation from the one-dimensional representation of $H_n^{aff}$ defined by the following,
\[
T_i \mapsto \zeta^{1/2}, \quad Y_j \mapsto \zeta^{\rho_j}.
\]
Here 
\[
\rho=(\rho_1, \cdots ,\rho_n)=\left(-\frac{n-1}{2},-\frac{n-3}{2}, \cdots ,\frac{n-1}{2}\right).
\]
\end{prop}
\begin{rem}
If $\zeta,\tau$ are generic, this polynomial representation is irreducible and $Y$-semisimple, namely the action of the commutative operators $Y_i \ (1 \le i \le n)$ are simultaneously diagonalizable. The simultaneous eigenvectors for $Y_i$, are the non-symmetric Macdonald polynomials. For more details, see \cite{Kasa}.
\end{rem}
\subsubsection{Category $\mathcal{O}$ of DAHA}
\indent \indent We introduce the $Y$-intertwining operators as follows
\[
\phi_i=T_i(Y_{\alpha_i}-1)+(\zeta^{1/2}-\zeta^{-1/2}).
\]
Let us define $\phi_w=\pi^m\phi_{s_{1_1}}\phi_{s_{i_2}} \cdots \phi_{s_{i_k}}$ if $w=\pi^ms_{i_1}s_{i_2} \cdots s_{i_k}$ is a reduced expression of $w \in W_n$. It does not depend on the choise of reduced expressions.
\begin{prop} \quad {} \\
(1) \ $\phi_wY_\eta=Y_{w\eta}\phi_w$ for any $w \in W_n, \eta \in P$. \\
(2) \ $\phi_i^2=\{\zeta^{1/2}(Y_{\alpha_i}-1)+(\zeta^{1/2}-\zeta^{-1/2}\}\{\zeta^{1/2}(Y_{-\alpha_i}-1)+(\zeta^{1/2}-\zeta^{-1/2}\}$. \\
(3) \ Put $\Phi_i=\{\zeta^{1/2}(Y_{\alpha_i}-1)+(\zeta^{1/2}-\zeta^{-1/2}\}^{-1}\phi_i$ and $\Phi_w=\pi^m\Phi_{s_{i_1}}\Phi_{s_{i_2}} \cdots \Phi_{s_{i_k}}$ if $w=\pi^ms_{i_1}s_{i_2} \cdots s_{i_k}$ is a reduced expression of $w \in W_n$. Then the operators $\{\Phi_w|w \in W_n\}$ satisfy the defining relation of extended affine Weyl group. 
\end{prop}

Let $\mathcal{H}_n$-mod be the category of finitely generated $\mathcal{H}_n$-modules. Put $\bb{C}[Y]$ the subalgebra of $\mathcal{H}_n$ generated by $Y_{\eta} \ (\eta \in P)$.
\begin{dfn}
The category $\mathcal{O}$ is the full subcategory of $\mathcal{H}_n$-mod consisting of modules which are locally finite with respect to $\bb{C}[Y]$. Here a module $M \in \mathcal{H}_n$-mod is locally finite with respect to $\bb{C}[Y]$ if $\bb{C}[Y]v$ is finite-dimensional for any $v \in M$.
\end{dfn}
Note that the polynomial representation of $\mathcal{H}_n$ belongs to $\mathcal{O}$. 
\quad {} \\
\indent Suppose $M \in \mathcal{O}$. Then $M$ has a generalized weight decomposition
\begin{eqnarray*}
M&=&\bigoplus_{\chi \in \mf{h}^*}M_{\chi}, \\
&{}&\text{where} \ M_{\chi}=\bigcup_{k \ge 1}\left\{v \in M|(Y_\eta-\zeta^{(\eta|\chi)})^kv=0  \ \text{for any} \  \eta \in P\right\}.
\end{eqnarray*}
Let us set $\Supp(M)=\{\chi \in \mf{h}^*|M_{\chi} \neq 0\}$.
\begin{dfn}
The category ${}^{\chi}\mathcal{O}$ is the full subcategory consisting of the modules $M$ such that $\Supp(M) \subset W_n \cdot \chi$.
\end{dfn}
Note that the polynomial representation of $\mathcal{H}_n$ belongs to ${}^{\rho}\mathcal{O}$.\\
\subsubsection{Specialized parameters}
\indent \indent In this paper, we will specialize the two parameters $\zeta$ and $\tau$
at 
\[
\zeta^\ell\tau^r=1 \quad (2 \le \ell \le n, \ 1 \le r, \ (\ell,r)=1).
\]
We assume that $\zeta,\tau$ are not roots of unity. Let $\mathcal{H}_n^{(\ell,r)}$ be the algebra $\mathcal{H}_n$ with the parameters specialized as above, and $V_n^{(\ell,r)}$ the polynomial representation of $\mathcal{H}_n^{(\ell,r)}$. Generally, the representation $V_n^{(\ell,r)}$ is not irreducible and not $Y$-semisimple in general.\\
\indent Let us define the full subcategories $\mathcal{O}_{(\ell,r)}$ and ${}^{\chi}\mathcal{O}_{(\ell,r)}$ of $\mathcal{H}_n^{(\ell,r)}$-mod similarly to the generic case.

\subsection{Kasatani's Conjecture on the polynomial representation of DAHA}
M. Kasatani constructed in \cite{Kasa} the subrepresentations of $V_n^{(\ell,r)}$ using ``wheels'' of variables with length $\ell$. This is called the ``multi-wheel condition''. In this section, we will recall his construction based on \cite{Kasa}. \\
\begin{dfn}
Let $Z_m^{(\ell,r)}$ be the subset of $\bb{K}^n$ contained in $(z_1, \cdots ,z_n) \in \bb{K}^n$ satisfying the following conditions;
\begin{quote}
there exist distinct indices
\[
i_{j,1}, \cdots ,i_{j,\ell} \in \{1, \cdots ,n\} \quad (1 \le j \le m)
\]
and positive integers
\[
s_{j,1}, \cdots ,s_{j,\ell} \in \bb{Z}_{\ge 0} \quad (1 \le j \le m)
\]
such that
\begin{eqnarray*}
&{}&z_{i_{j,a}}=\zeta\tau^{s_{j,a}}z_{i_{j,a+1}} \quad (1 \le j \le m, \ 1 \le i \le \ell), \\
&{}&\sum_{a=1}^{\ell}s_{j,a}=r \quad (1 \le j \le m), \\
&{}&i_{j,a+1}>i_{j,a} \ \text{if} \ s_{j,a}=0.
\end{eqnarray*}
\end{quote}
Let us define the ideals
\[
I_m^{(\ell,r)}=\{f \in \bb{K}[x_1^{\pm{1}}, \cdots ,x_n^{\pm{1}}]; f(z)=0 \ \text{for all} \  z \in Z_m^{(\ell,r)}\}.
\]
We call the defining relation of $I_m^{(\ell,r)}$ the multi-wheel condition.
\end{dfn}
We introduce Kasatani's result and conjecture.
\begin{thm}[{\cite[Theorem 6.3]{Kasa}}]
Let $N=\left[\frac{n}{\ell}\right]$. The sequence 
\[
0=I_0^{(\ell,r)} \subsetneq I_1^{(\ell,r)} \subsetneq  I_2^{(\ell,r)} \subsetneq \cdots I_N^{(\ell,r)} \subsetneq I_{N+1}^{(\ell,r)}=V_n^{(\ell,r)}.
\]
is an increasing sequence of subrepresentations of $V_m^{(\ell,r)}$.
\end{thm}
We call the following conjecture the Kasatani Conjecture on the polynomial representation of DAHA.
\begin{conj}[{\cite[Conjecture 6.4]{Kasa}}]
The above increasing sequence of subrepresentations of $V_n^{(\ell,r)}$ is a composition series, namely the quotient representations
\[
I_{a+1}^{(\ell,r)}/I_{a}^{\ell,r)} \quad (0 \le a \le N)
\]
are irreducible. Especially note that the number of composition factors of $V_n^{(\ell,r)}$ is equal to $N+1=\left[\frac{n}{\ell}\right]+1$.
\end{conj}
\newpage

\section{Two Degenerate Versions of DAHA and their category $\bb{O}$}
\subsection{Trigonometric Degeneration of DAHA}
\indent \indent In this section, we will recall two degenerate versions of DAHA and their category $\bb{O}$. They are called the ``trigonometric degeneration of DAHA'' and ``rational degeneration of DAHA''.
\subsubsection{Trigonometric degeneration of DAHA}
Let us define the trigonometric degeneration of DAHA, simply called degenerate DAHA.
\begin{dfn}
The degenerate DAHA $\bb{H}^{deg}_{n,h}$ of type $GL_n$ is an associative algebra over $\bb{C}$ generated by
\[
\pi^{\pm{1}}, s_0, s_1, \cdots s_{n-1}, \quad y^{deg}_\eta \ (\eta \in P \oplus \bb{Z}\delta)
\]
satisfying the following relations
\begin{eqnarray*}
\begin{array}{ll}
y^{deg}_{\delta}=1, & \\
y^{deg}_{\eta}+y^{deg}_{\xi}=y^{deg}_{\eta+\xi} & (\eta,\xi \in P), \\
\langle \pi^{\pm{1}}, s_0, s_1, \cdots ,s_{n-1} \rangle \cong \bb{C}W_n, & \\
s_iy^{deg}_\eta-y^{deg}_{s_i\eta}s_i=h\dfrac{y^{deg}_{s_i\eta}-y^{deg}_{\eta}}{y_{\alpha_i}} & (0 \le i \le n-1, \ \eta \in P), \\
{\pi}y^{deg}_{\eta}=y^{deg}_{\pi\eta}\pi, \\
\end{array}
\end{eqnarray*}
for $h \in \bb{C} \backslash \{0\}$.
\end{dfn}
Recall the another generators $X_1, \cdots ,X_n,s_1, \cdots ,s_n$ of the extended affine Weyl group $W_n$, namely
\[
X_1={\pi}s_{n-1}s_{n-2} \cdots s_1, \quad X_i=\pi^{i-1}X_1\pi^{-i+1}.
\]

\subsubsection{Standard modules and Category $\bb{O}^{deg}$ of degenerate DAHA}
\indent \indent We introduce the $Y$-intertwining operators as follows
\[
\phi_i^{deg}=s_iy_{\alpha_i}^{deg}+h.
\]
Let us define $\phi_w^{deg}=\pi^m\phi_{i_1}^{deg}\phi_{i_2}^{deg} \cdots \phi_{i_k}^{deg}$ if $w=\pi^ms_{i_1}s_{i_2} \cdots s_{i_k}$ is a reduced expression of $w \in W_n$. It does not depend on the choice of the reduced expressions.
\begin{prop} \quad {} \\
(1) \ $\phi_w^{deg}y_\eta^{deg}=y_{w\eta}^{deg}\phi_w^{deg}$ for any $w \in W_n, \eta \in P$. \\
(2) \ $(\phi_i^{deg})^2=(h-y^{deg}_{\alpha_i})(h+y^{deg}_{\alpha_i})$. \\
(3) \ Put $\Phi_i^{deg}=(h-y^{deg}_{\alpha_i})^{-1}\phi_i^{deg}$ and $\Phi_w^{deg}=\pi^m\Phi_{s_{i_1}}^{deg}\Phi_{s_{i_2}}^{deg} \cdots \Phi_{s_{i_k}}^{deg}$ for a reduced expression $w=\pi^ms_{i_1}s_{i_2} \cdots s_{i_k}$ of $w \in W_n$. Then the operators $\{\Phi_w^{deg}|w \in W_n\}$ satisfy the defining relation of extended affine Weyl group. 
\end{prop}

Let $\bb{H}^{deg}_{n,h}$-mod be the category of finitely generated $\bb{H}^{deg}_{n,h}$-modules. Let $\bb{C}[y^{deg}]$ be the subalgebra of $\bb{H}^{deg}_{n,h}$ generated by $y^{deg}_{\eta} \ (\eta \in P)$.
\begin{dfn}
The category $\bb{O}^{deg}_h$ is the full subcategory of $\bb{H}^{deg}_{n,h}$-mod consisting of modules which are locally finite with respect to $\bb{C}[y^{deg}]$. 
\end{dfn}
\quad {} \\
Suppose $M \in \bb{O}^{deg}_h$. Then $M$ has the generalized weight decomposition
\begin{eqnarray*}
M&=&\bigoplus_{\chi \in \mf{h}^*}M_{\chi}, \\
&{}&\text{where} \ M_{\chi}=\bigcup_{k \ge 1}\left\{v \in M|(y_\eta-(\eta|\chi))^kv=0  \ \text{for any} \ \eta \in P\right\}.
\end{eqnarray*}
Let us set $\Supp(M)=\{\chi \in \mf{h}^*|M_{\chi} \neq 0\}$.
\begin{dfn}
The category ${}^{\chi}\bb{O}^{deg}_h$ is the full subcategory consisting of the modules $M$ such that $\Supp(M) \subset W_n \cdot \chi$.
\end{dfn}
\quad {} \\
The degenerate DAHA $\bb{H}^{deg}_n$ has the polynomial representation on $\bb{C}[x_1^{\pm{1}}, \cdots ,x_n^{\pm{1}}]$. More generally, we can introduce the standard modules.
\begin{dfn}
Let $S^\lambda$ be the $\mf{S}_n$-module corresspond to a partition $\lambda$ of $n$. Then $S^\lambda$ becomes $\bb{C}\mf{S}_n \otimes \bb{C}[y^{deg}]$-module in which by the action of $y^{deg}$ given by
\[
y_i^{deg} \mapsto \sum_{j<i}s_{ji}-\frac{n-1}{2}.
\]
The standard module $\Delta_h^{deg}(\lambda)$ is the induced module of the $\bb{S}\mf{S}_n \otimes \bb{C}[y^{deg}]$-module $S^\lambda$ to $\bb{H}_{n,h}^{deg}$;
\[
\Delta_h^{deg}(\lambda)=\Ind_{\bb{C}\mf{S}_n \otimes \bb{C}[y^{deg}]}^{\bb{H}^{deg}_n}S^\lambda.
\]
Especially, the standard module $\Delta_h^{deg}(\triv)$ is isomorphic to the polynomial representation of $\bb{H}_{n,h}^{deg}$. Namely, the polynomial representation of $\bb{H}^{deg}_{n,h}$ is induced module from the one-dimensional $\bb{C}\mf{S}_n \otimes \bb{C}[y^{deg}]$-module:
\[
s_i \mapsto 1 \ (1 \le i \le n-1), \quad y^{deg}_j \mapsto \rho_j \ (1 \le j \le n).
\]
Recall that $\rho=(\rho_1, \cdots ,\rho_n)=\left(-\frac{n-1}{2},-\frac{n-3}{2}, \cdots ,\frac{n-1}{2}\right)$.
\end{dfn}
Note that the polynomial representation $\Delta_h^{deg}(\triv)$ belongs to the category ${}^{\rho}\bb{O}^{deg}_{h}$. 
\begin{prop}[\cite{Su}]
Each standard module $\Delta^{deg}_h(\lambda)$ has a unique simple quotient denoted by $L^{deg}_h(\lambda)$.
\end{prop}
\subsection{Rational Degeneration of DAHA}
\subsubsection{Rational degeneration of DAHA}
Let us define the rational degeneration of DAHA, simply called rational DAHA.
\begin{dfn}
The degenerate DAHA $\bb{H}^{rat}_{n,h}$ of type $GL_n$ is the associative algebra over $\bb{C}$ generated by
\[
x_{\eta^\vee} \ (\eta^\vee \in P^{\vee}), \quad s_1, \cdots ,s_{n-1}, \quad y^{rat}_\eta \ (\eta \in P)
\]
satisfying the following relations
\begin{eqnarray*}
&{}&\begin{array}{ll}
y^{rat}_{\eta}+y^{rat}_{\xi}=y^{rat}_{\eta+\xi} & (\eta,\xi \in P) \\
x_{\eta^\vee}+x_{\xi^\vee}=x_{\eta^\vee+\xi^\vee} & (\eta^\vee,\xi^\vee \in P^\vee) \\
\langle s_1, \cdots ,s_{n-1} \rangle \cong \bb{C}\mf{S}_n, & \\
wx_{\eta^\vee}=x_{w\eta^\vee}w & (w \in \mf{S}_n), \\
wy_\eta^{rat}=y_{w\eta}^{rat}w, & 
\end{array} 
\quad {} \\
&{}&[x_i,y_j^{rat}]=
\left\{
\begin{array}{ll}
hs_{ij} & (\text{if} \ i \neq j) \\
\displaystyle 1-h\sum_{k \neq i}{s_{ik}} & (\text{if} \ i=j)
\end{array}
\right.,
\end{eqnarray*}
for $h \in \bb{C} \backslash \{0\}$, where $x_i=x_{\varepsilon_i^{\vee}},y^{rat}_i=y^{rat}_{\varepsilon_i}$.

\end{dfn}

\subsubsection{Standard modules and Category $\bb{O}^{rat}$ of rational DAHA}
Let $\bb{H}^{rat}_{n,h}$-mod be the category of finitely generated $\bb{H}^{rat}_{n,h}$-modules. Let $\bb{C}[y^{rat}]$ be the subalgebra of $\bb{H}^{rat}_{n,h}$ generated by $y^{rat}_{\eta} \ (\eta \in P)$.
\begin{dfn}
The category $\bb{O}^{rat}_h$ is the full subcategory of $\bb{H}^{rat}_{n,h}$-mod consisting of modules which are locally nilpotent with respect to $y^{rat}$. Here a module $M \in \bb{H}^{rat}_{n,h}$-mod is locally nilpotent with respect to $y^{rat}$ if for any $v \in M$ there exists $N>0$ such thet $(y_i^{rat})^Nv=0 \ (1 \le i \le n)$.
\end{dfn}
The rational DAHA $\bb{H}^{rat}_n$ has the polynomial representation on $\bb{C}[x_1, \cdots ,x_n]$. More generally, we can introduce the standard modules.
\begin{dfn}
Let $S^\lambda$ be the $\mf{S}_n$-module corresspond to a partition $\lambda$ of $n$. Then $S^\lambda$ becomes a $\bb{C}\mf{S}_n \otimes \bb{C}[y^{rat}]$-module by
\[
y_i^{rat}S^{\lambda}=0.
\]
The standard module $\Delta_h^{rat}(\lambda)$ is the $\bb{H}_{n,h}^{rat}$-module induced by the $\bb{C}\mf{S}_n \otimes \bb{C}[y^{rat}]$-module;
\[
\Delta_h^{rat}(\lambda)=\Ind_{\bb{C}\mf{S}_n \otimes \bb{C}[y^{rat}]}^{\bb{H}^{rat}_n}S^\lambda.
\]
The standard module $\Delta_h^{rat}(\triv)$ is isomorphic to the polynomial representation of $\bb{H}_{n,h}^{rat}$. 
\end{dfn}
Note that the polynomial representation $\Delta_h^{rat}(\triv)$ belongs to the category $\bb{O}^{rat}_{h}$. \\
\begin{prop} \quad {} \\
(1) \ The standard modules $\Delta^{rat}(\lambda)$ have a unique simple quotient denoted by $L^{rat}(\lambda)$. \\
(2) \ The set $\{L^{rat}(\lambda)\}$ is a complete collection of irreducible representations of $\bb{H}^{rat}_h$. \\
\end{prop}
\subsubsection{Embedding to degenerate DAHA}
The rational DAHA can be embedded to the degenerate DAHA proved by T. Suzuki.
\begin{prop}[\cite{Su}]
The following homomorphism from $\bb{H}_{n,h}^{rat}$ to $\bb{H}_{n,h}^{deg}$ is an embedding;
\begin{eqnarray*}
s_i &\mapsto& s_i, \\
x_j^{\vee} &\mapsto& X_j=\pi^{j-1}X_1\pi^{-j+1} \ \text{where} \  X_1={\pi}s_{n-1} \cdots s_1, \\
y_j^{rat} &\mapsto& X_j^{-1}\left(y_j^{deg}-\sum_{1 \le k<j}s_{kj}+\frac{n-1}{2}\right).
\end{eqnarray*}
\end{prop}
\newpage

\section{Relationships between Category $\mathcal{O}$ and $v$-Schur algebras}
\indent \indent In this section, we will explain the relationships of some categories.
\subsection{Relationship of $\mathcal{O}$ and $\bb{O}^{deg}$}
\indent \indent The categories $\mathcal{O}$ and $\bb{O}^{deg}$ are not equivalent. But there exists an equivalence of categories between their full subcategories when the parameters are special. The following theorem proved by Varagnolo-Vasserot and Lusztig.
\begin{thm}[\cite{VV2},\cite{Lus}]
If $\chi \in \mf{h}$ satisfies that
\[
(\eta|\chi) \in \bb{Z} \ \text{and} \  (\delta|\chi) \in \bb{Z} \ \text{for any} \ \eta \in P,
\]
then the categories ${}^{\chi}\mathcal{O}_{(\ell,r)}$ of the specialized DAHA $\mathcal{H}_{n}^{(\ell,r)}$ and ${}^{\chi}\bb{O}_{r/\ell}^{deg}$ of the degenerate DAHA $\bb{H}_{n,r/\ell}^{deg}$ are equivalent.
\end{thm} 
We sketch the proof of this theorem. \\
\indent First, for any finite subset $E \subset W_n \cdot \chi$, put 
\[
\langle E \rangle=\bigcap_{\xi \in E}\{Y_\eta-\zeta^{(\eta|\xi)}|\eta \in P\}.
\]
We consider the $\langle{E}\rangle$-adic completion $S[Y]$ of $\bb{C}[Y]$. Let $S(Y)$ be the quotient field of $S[Y]$. Similarly, put
\[
\langle{E}\rangle^{deg}=\bigcap_{\xi \in E}\{y_\eta^{deg}-{(\eta|\xi)}|\eta \in P\},
\]
and consider the $\langle{E}\rangle^{deg}$-adic completion $S[y^{deg}]$ of $\bb{C}[y^{deg}]$. Let $S(y^{deg})$ be the quotient field of $S[y^{deg}]$. Then there exists an isomorphism $S[Y] \cong S[y^{deg}]$, and this isomorphism can be extended $j:S(Y) \to S(y^{deg})$. Therefore there exists a morphism
\[
\widetilde{j}:S(Y) \otimes_{\bb{C}[Y]}\mathcal{H}_n^{(\ell,r)} \to S(y^{deg}) \otimes_{\bb{C}[y^{deg}]} \bb{H}_{n,r/\ell}^{deg}
\]
such that 
\[
\sum_{w \in W_n}c_w\Phi_w \mapsto \sum_{w \in W_n}j(c_w)\Phi_w^{deg}.
\]
Then this morphism induce the algebra isomorphism of $S[Y] \otimes_{\bb{C}[Y]}\mathcal{H}_n^{(\ell,r)} \to S[y^{deg}] \otimes_{\bb{C}[y^{deg}]} \bb{H}_{n,r/\ell}^{deg}$. \\
\indent The category of "smooth" modules of $S[Y] \otimes_{\bb{C}[Y]}\mathcal{H}_n^{(\ell,r)}$ is equivalent to the category ${}^{\chi}\mathcal{O}_{(\ell,r)}$. The category of "smooth" modules of $S[y^{deg}] \otimes_{\bb{C}[y^{deg}]}\bb{H}_{n,r/\ell}^{deg}$ is equivalent to the category ${}^{\chi}\bb{O}_{r/\ell}$. Thus the categories ${}^{\chi}\mathcal{O}_{(\ell,r)}$ and ${}^{\chi}\bb{O}_{r/\ell}$ are equivalent. For more details, see \cite{VV2},\cite{Lus}.

\subsection{Embedding of $\bb{O}^{rat}$ to $\bb{O}^{deg}$}
By the embedding $\bb{H}_{n,h}^{rat} \hookrightarrow \bb{H}_{n,h}^{deg}$, we can define the induction functor $\bb{H}_{n,h}^{rat}$-mod to $\bb{H}_{n,h}^{deg}$-mod. The following theorem states that this functor is fully faithful and exact.

\begin{thm}[\cite{Su}] \quad {} \\
(1) \ The functor 
\[
\bb{O}^{rat}_{h} \to \bb{O}_{h}^{deg}; M \mapsto \bb{H}_{n,h}^{deg} \otimes_{\bb{H}^{rat}_{n,h}}M
\]
is fully faithful and exact.\\
(2) \ The above functor sends the standard modules to standard modules, namely 
\[
\Delta^{deg}(\lambda)=\bb{H}_{n,h}^{deg} \otimes_{\bb{H}^{rat}_{n,h}}\Delta^{rat}(\lambda).
\]
Especially
\[
\Delta^{deg}(\triv)=\bb{H}_{n,h}^{deg} \otimes_{\bb{H}^{rat}_{n,h}}\Delta^{rat}(\triv).
\]
(3) \  The above functor sends the simple module $L^{rat}(\lambda)$ to $L^{deg}(\lambda)$, namely
\[
L^{deg}(\lambda)=\bb{H}_{n,h}^{deg} \otimes_{\bb{H}^{rat}_{n,h}}L^{rat}(\lambda).
\]
Thus we obtain $[\Delta^{deg}(\lambda):L^{deg}(\mu)]=[\Delta^{rat}(\lambda):L^{rat}(\mu)]$.

\end{thm}
\subsection{$v$-Schur algebra $\bb{S}(n)$}
\indent \indent The $v$-Schur algebra was introduced by Dipper and James \cite{DJ}. Their definition used the Hecke algebra of type $A$. On the other hand, Beilinson-Lusztig-Macpherson (\cite{BLM}) constructed the $v$-Schur algebra by using a geometry of flag varieties. Their construction is related to the quantum enveloping algebra and the Schur-Weyl duality.\\
\indent Let $H_n$ be the Iwahori-Hecke algebra of the symmetric group $\mf{S}_n$ with a parameter $v$. For a composition $\mu=(\mu_1, \cdots ,\mu_n)$ of $n$, let us denote by $\mf{S}_{\mu}$ the Young subgroup $\mf{S}_{\mu_1} \times \cdots \times \mf{S}_{\mu_n}$. Put $\displaystyle m_{\mu}=\sum_{w \in \mf{S}_\mu}T_w \in H_n$. Then the left $H_n$-module 
\[
M=\bigoplus_{\mu}H_nm_{\mu}
\]
is called the permutation module. The $v$-Schur algebra $\bb{S}(n)$ is the endomorphism ring of $M$;
\[
\bb{S}(n) \cong \End_{H_n}(M)^{\text{op}}.
\]
Moreover the $H_n$-module $M$ naturally becomes an $\bb{S}(n)$-module. Thus we can define the following functor;
\[
\mathcal{S}:\bb{S}(n)\text{-mod} \to H_n\text{-mod}; \ \text{given by} \ N \mapsto M \otimes_{\bb{S}(n)} N.
\]
If $v$ is not a root of unity, we can construct a complete collection of irreducible $\bb{S}(n)$-modules $\{W^\lambda|\lambda \vdash n\}$ indexed by the partitions of $n$. The irreducible module $W^\lambda$, called the Weyl module, satisfies
\[
S^{\lambda}=M \otimes_{\bb{S}(n)} W^\lambda.
\]
Here, the $H_n$-modules $\{S^\lambda|\lambda \vdash n\}$ is a complete collection of irreducible $H_n$-modules when $v$ is not a root of unity. \\ 
\indent It is known that the $v$-Schur algebra has a cellular algebra structure, a notion  introduced by Graham-Lehrer \cite{GL}. Therefore if $v$ is a root of unity,  $L^\lambda=W^\lambda/\rad{W^\lambda}$ are irreducible unless $L^\lambda=0$. Moreover, $L^\lambda \neq 0$ for any partititon $\lambda$ of $n$. Thus
\[
\{L^\lambda|\lambda \vdash n\}
\]
 is a complete set of irreducible representation of $\bb{S}(n)$. \\
\indent On the other hand, we consider the quantum enveloping algebra $U_v(\mf{gl}_n)$ and its vector representation $\bb{C}^n$. Let $(\bb{C}^n)^{\otimes n}$ be the quantum tensor product representation;
\[
U_v(\mf{gl}_n) \to \End((\bb{C}^n)^{\otimes n}).
\]
Then the $v$-Schur algebra $\bb{S}(n)$ is isomorphic to the image of this homomorphism. \\
\indent The $U_q(\widehat{\mf{gl}_n})$-module associated with the $\bb{S}(n)$-module $W^{(1^n)}$ by the above quotient map is the determinant representation. Note that if $v$ is not a root of unity, the above functor $\mathcal{S}$ sends the irreducible $\bb{S}(n)$-module $W^{(1^n)}$ to the irreducible $H_n$-module $S^{(1^n)}$ called the signature representation.

\subsection{Equivalence of $\bb{S}(n)$-mod and $\bb{O}^{rat}$}
\indent \indent R. Rouquier proved in \cite{Rou2} that there exists an equivalence of categories between $\bb{O}^{rat}$ and $\bb{S}(n)$-mod. 
\begin{thm}
Let us consider the categories $\bb{O}^{rat}_h$ and $\bb{S}(n)$-mod at $v=\sqrt[\ell]{1}$. When $h \notin \dfrac{1}{2}+\bb{Z}$, there exists an equivalence of categories
\[
\Psi^{Schur} \ : \ \bb{S}(n)\text{-mod} \ \stackrel{\sim}{\longrightarrow} \ \bb{O}^{rat}_h
\]
such that if $h>0$ the standard modules $\Delta_h^{rat}(\lambda)$ are sent to the Weyl module $W^{\lambda'}$. Her $\lambda'$ is the conjugate partition of $\lambda$.
\end{thm}

\subsection{Summary on the case of polynomial representation}
 \indent \indent We defined three functors in preceding sections:
\begin{eqnarray*}
{}^{\rho}\Psi^{deg}&:& {}^{\rho}\bb{O}_{r/\ell}^{deg} \stackrel{\sim}{\longrightarrow} {}^{\rho}\mathcal{O}_{r/\ell}, \\
\Psi^{rat}&:&\bb{O}_{r/\ell}^{rat} \longrightarrow \bb{O}_{r/\ell}^{deg}, \\
\Psi^{Schur}&:&\bb{S}(n)\text{-mod} \stackrel{\sim}{\longrightarrow} \bb{O}_{r/\ell}^{rat}.
\end{eqnarray*}
The two functors ${}^{\rho}\Psi^{deg}$ and $\Psi^{Schur}$ are equivalences of categories, and the functor $\Psi^{rat}$ is fully faithful and exact. The representation 
\[
{}^{\rho}\Psi^{deg}(\Psi^{rat} \circ \Psi^{Schur}(W^{(1^n)}))
\] 
is isomorphic to the polynomial representation of $\mathcal{H}_{n}^{(\ell,r)}$. \\
 \indent We can reduce the determination of the composition factors of the polynomial representation of $\mathcal{H}_n^{(\ell,r)}$ to the one of the composition factors of $\bb{S}(n)$-module $W^{(1^n)}$ by using the above correspondence.
\newpage
\section{Global and Crystal basis of Fock Space}
\subsection{Quantum enveloping algebra $U_q(\widehat{\mf{sl}_\ell})$}
Let us recall the quantum enveloping algebra $U_q(\widehat{\mf{sl}_\ell})$. Set $[n]=\dfrac{q^n-q^{-n}}{q-q^{-1}}$. Let the matrix $A=(a_{ij})_{0 \le i \le \ell-1}$ be the Cartan matrix of type $A_{\ell-1}^{(1)}$. Namely, if $\ell \ge 3$, 
\[
a_{ij}=\left\{
\begin{array}{ll}
2 & i=j \\
-1 & i \equiv j \pm{1} (\mod{\ell})\\
0 & \text{otherwise} 
\end{array}
\right.
\]
and if $\ell=2$
\[
A=\left(
\begin{array}{cc}
2 & -2 \\
-2 & 2
\end{array}
\right).
\]
\begin{dfn}
The quantum enveloping algebra $U_q(\widehat{\mf{sl}_\ell})$ is generated by 
\[
E_i, F_i, K_i \quad (0 \le i \le \ell-1),
\]
satisfying the following relations
\begin{eqnarray*}
&{}&K_iK_j=K_jK_i, \\
&{}&K_iE_j=q^{a_{ij}}E_jK_i, \\
&{}&K_iF_j=q^{-a_{ij}}F_jK_i, \\
&{}&E_iF_j-F_jE_i=\delta_{i,j}\frac{K_i-K_i^{-1}}{q-q^{-1}}, \\
&{}&E_iE_j=E_jE_i \quad (\text{if} \ i \neq j \pm{1}), \\
&{}&F_iF_j=F_jF_i \quad (\text{if} \ i \neq j \pm{1}), 
\end{eqnarray*}
and the $q$-Serre relations,
\begin{eqnarray*}
\text{if} \ \ell \ge 3, &{}& \\
&{}&E_i^2E_{i\pm{1}}-(q+q^{-1})E_iE_{i\pm{1}}E_i+E_{i\pm{1}}E_i^2=0, \\
&{}&F_i^2F_{i\pm{1}}-(q+q^{-1})F_iF_{i\pm{1}}F_i+F_{i\pm{1}}F_i^2=0, \\
\text{if} \ \ell=2, &{}& \\
&{}&E_i^3E_{i\pm{1}}-[3]E_i^2E_{i\pm{1}}E_i+[3]E_iE_{i\pm{1}}E_i^2-E_{i\pm{1}}E_i^2=0, \\
&{}&F_i^3F_{i\pm{1}}-[3]F_i^2F_{i\pm{1}}F_i+[3]F_iF_{i\pm{1}}F_i^2-F_{i\pm{1}}E_i^2=0,
\end{eqnarray*}
The indices in the above relations are to be read modulo $\ell$.
\end{dfn}
\indent $U_q(\widehat{\mf{sl}_\ell})$ is a Hopf algebra with a coproduct given by
\begin{eqnarray*}
\Delta^-(E_i)&=&1 \otimes E_i+E_i \otimes K_i^{-1}, \\
\Delta^-(F_i)&=&F_i \otimes 1+K_i \otimes F_i, \\
\Delta^-(K_i)&=&K_i \otimes K_i.
\end{eqnarray*}
There exists another coproduct $\Delta^+$ on $U_q(\widehat{\mf{sl}_n})$ given by
\begin{eqnarray*}
\Delta^+(E_i)&=&E_i \otimes 1+K_i \otimes E_i, \\
\Delta^+(F_i)&=&F_i \otimes K_i^{-1}+1 \otimes F_i, \\
\Delta^+(K_i)&=&K_i \otimes K_i.
\end{eqnarray*}
When we consider the lower and upper golobal basis, we will use these two coproducts.\\
\subsection{Two realizations of Fock space}
In this subsection, we will describe two realizations of the Fock space of $U_q(\widehat{\mf{sl}_\ell})$, the ``Hayashi realization'' and semi-infinite wedge space. 
\subsubsection{Hayashi realization}
\indent \indent Let us recall some notations and definitions. \\
\indent A partition $\lambda=(\lambda_1,\lambda_2, \cdots )$ is a non-increasing sequence of natural numbers. The corresponding Young diagram is a collection of rows of square boxes which are left justified and $\lambda_i$ boxes in the $i$-th row. Let $\mathcal{P}$ be the set of all partitions. 
\begin{dfn} \quad {} \\
(1) \ The content of box $x \in \lambda$ is defined by 
\[
c(x)=\col(x)-\row(x).
\]
Assume that we are given a positive number $\ell$. The $\ell$-residue of a box $x$ is defined by $\res_\ell(x)=c(x)$ modulo $\ell$.\\
(2) \ If a partiton $\mu$ is obtained by removing a box $x$ from a partition $\lambda$, the box $x$ is called a removable box in $\lambda$. Conversely, if a partition $\lambda$ is obtained by adding a box $x$ to a partition $\mu$, the box $x$ is call an addable box in $\mu$. If a removable [resp. addable] box $x$ has residue $i$, we call the box $x$ an $i$-removable [resp. $i$-addable] box. \\ 
\end{dfn}
\[
\unitlength 0.1in
\begin{picture}(62.00,35.70)(6.00,-36.00)
%
\special{pn 8}%
\special{pa 3000 400}%
\special{pa 3400 400}%
\special{pa 3400 800}%
\special{pa 3000 800}%
\special{pa 3000 400}%
\special{fp}%
%
\special{pn 8}%
\special{pa 3400 400}%
\special{pa 3800 400}%
\special{pa 3800 800}%
\special{pa 3400 800}%
\special{pa 3400 400}%
\special{fp}%
%
\special{pn 8}%
\special{pa 3800 400}%
\special{pa 4200 400}%
\special{pa 4200 800}%
\special{pa 3800 800}%
\special{pa 3800 400}%
\special{fp}%
%
\special{pn 8}%
\special{pa 4200 400}%
\special{pa 4600 400}%
\special{pa 4600 800}%
\special{pa 4200 800}%
\special{pa 4200 400}%
\special{fp}%
%
\special{pn 8}%
\special{pa 3000 800}%
\special{pa 3400 800}%
\special{pa 3400 1200}%
\special{pa 3000 1200}%
\special{pa 3000 800}%
\special{fp}%
%
\special{pn 8}%
\special{pa 3400 800}%
\special{pa 3800 800}%
\special{pa 3800 1200}%
\special{pa 3400 1200}%
\special{pa 3400 800}%
\special{fp}%
%
\special{pn 8}%
\special{pa 3800 800}%
\special{pa 4200 800}%
\special{pa 4200 1200}%
\special{pa 3800 1200}%
\special{pa 3800 800}%
\special{fp}%
%
\special{pn 8}%
\special{pa 3000 1200}%
\special{pa 3400 1200}%
\special{pa 3400 1600}%
\special{pa 3000 1600}%
\special{pa 3000 1200}%
\special{fp}%
\put(32.0000,-6.0000){\makebox(0,0){$0$}}%
\put(36.0000,-10.0000){\makebox(0,0){$0$}}%
\put(36.0000,-6.0000){\makebox(0,0){$1$}}%
\put(40.0000,-6.0000){\makebox(0,0){$2$}}%
\put(43.9000,-6.0000){\makebox(0,0){$3$}}%
\put(40.0000,-10.0000){\makebox(0,0){$1$}}%
\put(32.0000,-10.0000){\makebox(0,0){$-1$}}%
\put(32.0000,-14.0000){\makebox(0,0){$-2$}}%
%
\special{pn 8}%
\special{pa 3400 1200}%
\special{pa 3800 1200}%
\special{pa 3800 1600}%
\special{pa 3400 1600}%
\special{pa 3400 1200}%
\special{fp}%
\put(36.0000,-14.0000){\makebox(0,0){$-1$}}%
%
\special{pn 8}%
\special{pa 3400 1600}%
\special{pa 3800 1600}%
\special{pa 3800 2000}%
\special{pa 3400 2000}%
\special{pa 3400 1600}%
\special{fp}%
\put(36.0000,-18.0000){\makebox(0,0){$-2$}}%
%
\special{pn 8}%
\special{pa 3000 1600}%
\special{pa 3400 1600}%
\special{pa 3400 2000}%
\special{pa 3000 2000}%
\special{pa 3000 1600}%
\special{fp}%
\put(32.0000,-18.0000){\makebox(0,0){$-3$}}%
%
\special{pn 8}%
\special{pa 3000 2000}%
\special{pa 3400 2000}%
\special{pa 3400 2400}%
\special{pa 3000 2400}%
\special{pa 3000 2000}%
\special{fp}%
\put(32.0000,-22.0000){\makebox(0,0){$-4$}}%
%
\special{pn 8}%
\special{pa 5200 400}%
\special{pa 5600 400}%
\special{pa 5600 800}%
\special{pa 5200 800}%
\special{pa 5200 400}%
\special{fp}%
%
\special{pn 8}%
\special{pa 5600 400}%
\special{pa 6000 400}%
\special{pa 6000 800}%
\special{pa 5600 800}%
\special{pa 5600 400}%
\special{fp}%
%
\special{pn 8}%
\special{pa 6000 400}%
\special{pa 6400 400}%
\special{pa 6400 800}%
\special{pa 6000 800}%
\special{pa 6000 400}%
\special{fp}%
%
\special{pn 8}%
\special{pa 6400 400}%
\special{pa 6800 400}%
\special{pa 6800 800}%
\special{pa 6400 800}%
\special{pa 6400 400}%
\special{fp}%
%
\special{pn 8}%
\special{pa 5200 800}%
\special{pa 5600 800}%
\special{pa 5600 1200}%
\special{pa 5200 1200}%
\special{pa 5200 800}%
\special{fp}%
%
\special{pn 8}%
\special{pa 5600 800}%
\special{pa 6000 800}%
\special{pa 6000 1200}%
\special{pa 5600 1200}%
\special{pa 5600 800}%
\special{fp}%
%
\special{pn 8}%
\special{pa 6000 800}%
\special{pa 6400 800}%
\special{pa 6400 1200}%
\special{pa 6000 1200}%
\special{pa 6000 800}%
\special{fp}%
%
\special{pn 8}%
\special{pa 5200 1200}%
\special{pa 5600 1200}%
\special{pa 5600 1600}%
\special{pa 5200 1600}%
\special{pa 5200 1200}%
\special{fp}%
\put(54.0000,-6.0000){\makebox(0,0){$0$}}%
\put(58.0000,-10.0000){\makebox(0,0){$0$}}%
\put(58.0000,-6.0000){\makebox(0,0){$1$}}%
\put(62.0000,-6.0000){\makebox(0,0){$2$}}%
\put(65.9000,-6.0000){\makebox(0,0){$0$}}%
\put(62.0000,-10.0000){\makebox(0,0){$1$}}%
\put(54.0000,-10.0000){\makebox(0,0){$2$}}%
\put(54.0000,-14.0000){\makebox(0,0){$1$}}%
%
\special{pn 8}%
\special{pa 5600 1200}%
\special{pa 6000 1200}%
\special{pa 6000 1600}%
\special{pa 5600 1600}%
\special{pa 5600 1200}%
\special{fp}%
\put(58.0000,-14.0000){\makebox(0,0){$2$}}%
%
\special{pn 8}%
\special{pa 5600 1600}%
\special{pa 6000 1600}%
\special{pa 6000 2000}%
\special{pa 5600 2000}%
\special{pa 5600 1600}%
\special{fp}%
\put(58.0000,-18.0000){\makebox(0,0){$1$}}%
%
\special{pn 8}%
\special{pa 5200 1600}%
\special{pa 5600 1600}%
\special{pa 5600 2000}%
\special{pa 5200 2000}%
\special{pa 5200 1600}%
\special{fp}%
\put(54.0000,-18.0000){\makebox(0,0){$0$}}%
%
\special{pn 8}%
\special{pa 5200 2000}%
\special{pa 5600 2000}%
\special{pa 5600 2400}%
\special{pa 5200 2400}%
\special{pa 5200 2000}%
\special{fp}%
\put(54.0000,-22.0000){\makebox(0,0){$2$}}%
\put(6.0000,-2.0000){\makebox(0,0)[lb]{$\lambda=(4,3,2,2,1)$}}%
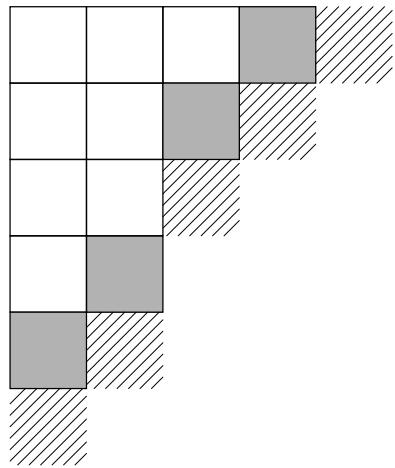
\put(30.0000,-26.0000){\makebox(0,0)[lt]{content of $\lambda$}}%
\put(51.9000,-26.0000){\makebox(0,0)[lt]{$3$-residue of $\lambda$}}%
%
\special{pn 8}%
\special{pa 600 400}%
\special{pa 1000 400}%
\special{pa 1000 800}%
\special{pa 600 800}%
\special{pa 600 400}%
\special{fp}%
%
\special{pn 8}%
\special{pa 1000 400}%
\special{pa 1400 400}%
\special{pa 1400 800}%
\special{pa 1000 800}%
\special{pa 1000 400}%
\special{fp}%
%
\special{pn 8}%
\special{pa 1400 400}%
\special{pa 1800 400}%
\special{pa 1800 800}%
\special{pa 1400 800}%
\special{pa 1400 400}%
\special{fp}%
%
\special{pn 8}%
\special{sh 0.300}%
\special{pa 1800 400}%
\special{pa 2200 400}%
\special{pa 2200 800}%
\special{pa 1800 800}%
\special{pa 1800 400}%
\special{fp}%
%
\special{pn 8}%
\special{pa 600 800}%
\special{pa 1000 800}%
\special{pa 1000 1200}%
\special{pa 600 1200}%
\special{pa 600 800}%
\special{fp}%
%
\special{pn 8}%
\special{pa 1000 800}%
\special{pa 1400 800}%
\special{pa 1400 1200}%
\special{pa 1000 1200}%
\special{pa 1000 800}%
\special{fp}%
%
\special{pn 8}%
\special{sh 0.300}%
\special{pa 1400 800}%
\special{pa 1800 800}%
\special{pa 1800 1200}%
\special{pa 1400 1200}%
\special{pa 1400 800}%
\special{fp}%
%
\special{pn 8}%
\special{pa 600 1200}%
\special{pa 1000 1200}%
\special{pa 1000 1600}%
\special{pa 600 1600}%
\special{pa 600 1200}%
\special{fp}%
%
\special{pn 8}%
\special{sh 0.300}%
\special{pa 1000 1600}%
\special{pa 1400 1600}%
\special{pa 1400 2000}%
\special{pa 1000 2000}%
\special{pa 1000 1600}%
\special{fp}%
%
\special{pn 8}%
\special{pa 600 1600}%
\special{pa 1000 1600}%
\special{pa 1000 2000}%
\special{pa 600 2000}%
\special{pa 600 1600}%
\special{fp}%
%
\special{pn 8}%
\special{sh 0.300}%
\special{pa 600 2000}%
\special{pa 1000 2000}%
\special{pa 1000 2400}%
\special{pa 600 2400}%
\special{pa 600 2000}%
\special{fp}%
%
\special{pn 8}%
\special{pa 1000 1200}%
\special{pa 1400 1200}%
\special{pa 1400 1600}%
\special{pa 1000 1600}%
\special{pa 1000 1200}%
\special{fp}%
%
\special{pn 8}%
\special{pa 1000 3000}%
\special{pa 1200 3000}%
\special{pa 1200 3200}%
\special{pa 1000 3200}%
\special{pa 1000 3000}%
\special{fp}%
%
\special{pn 4}%
\special{pa 1800 1220}%
\special{pa 1420 1600}%
\special{fp}%
\special{pa 1760 1200}%
\special{pa 1400 1560}%
\special{fp}%
\special{pa 1700 1200}%
\special{pa 1400 1500}%
\special{fp}%
\special{pa 1640 1200}%
\special{pa 1400 1440}%
\special{fp}%
\special{pa 1580 1200}%
\special{pa 1400 1380}%
\special{fp}%
\special{pa 1520 1200}%
\special{pa 1400 1320}%
\special{fp}%
\special{pa 1460 1200}%
\special{pa 1400 1260}%
\special{fp}%
\special{pa 1800 1280}%
\special{pa 1480 1600}%
\special{fp}%
\special{pa 1800 1340}%
\special{pa 1540 1600}%
\special{fp}%
\special{pa 1800 1400}%
\special{pa 1600 1600}%
\special{fp}%
\special{pa 1800 1460}%
\special{pa 1660 1600}%
\special{fp}%
\special{pa 1800 1520}%
\special{pa 1720 1600}%
\special{fp}%
%
\special{pn 4}%
\special{pa 2200 820}%
\special{pa 1820 1200}%
\special{fp}%
\special{pa 2160 800}%
\special{pa 1800 1160}%
\special{fp}%
\special{pa 2100 800}%
\special{pa 1800 1100}%
\special{fp}%
\special{pa 2040 800}%
\special{pa 1800 1040}%
\special{fp}%
\special{pa 1980 800}%
\special{pa 1800 980}%
\special{fp}%
\special{pa 1920 800}%
\special{pa 1800 920}%
\special{fp}%
\special{pa 1860 800}%
\special{pa 1800 860}%
\special{fp}%
\special{pa 2200 880}%
\special{pa 1880 1200}%
\special{fp}%
\special{pa 2200 940}%
\special{pa 1940 1200}%
\special{fp}%
\special{pa 2200 1000}%
\special{pa 2000 1200}%
\special{fp}%
\special{pa 2200 1060}%
\special{pa 2060 1200}%
\special{fp}%
\special{pa 2200 1120}%
\special{pa 2120 1200}%
\special{fp}%
%
\special{pn 4}%
\special{pa 2600 420}%
\special{pa 2220 800}%
\special{fp}%
\special{pa 2560 400}%
\special{pa 2200 760}%
\special{fp}%
\special{pa 2500 400}%
\special{pa 2200 700}%
\special{fp}%
\special{pa 2440 400}%
\special{pa 2200 640}%
\special{fp}%
\special{pa 2380 400}%
\special{pa 2200 580}%
\special{fp}%
\special{pa 2320 400}%
\special{pa 2200 520}%
\special{fp}%
\special{pa 2260 400}%
\special{pa 2200 460}%
\special{fp}%
\special{pa 2600 480}%
\special{pa 2280 800}%
\special{fp}%
\special{pa 2600 540}%
\special{pa 2340 800}%
\special{fp}%
\special{pa 2600 600}%
\special{pa 2400 800}%
\special{fp}%
\special{pa 2600 660}%
\special{pa 2460 800}%
\special{fp}%
\special{pa 2600 720}%
\special{pa 2520 800}%
\special{fp}%
%
\special{pn 4}%
\special{pa 1260 2000}%
\special{pa 1000 2260}%
\special{fp}%
\special{pa 1320 2000}%
\special{pa 1000 2320}%
\special{fp}%
\special{pa 1380 2000}%
\special{pa 1000 2380}%
\special{fp}%
\special{pa 1400 2040}%
\special{pa 1040 2400}%
\special{fp}%
\special{pa 1400 2100}%
\special{pa 1100 2400}%
\special{fp}%
\special{pa 1400 2160}%
\special{pa 1160 2400}%
\special{fp}%
\special{pa 1400 2220}%
\special{pa 1220 2400}%
\special{fp}%
\special{pa 1400 2280}%
\special{pa 1280 2400}%
\special{fp}%
\special{pa 1400 2340}%
\special{pa 1340 2400}%
\special{fp}%
\special{pa 1200 2000}%
\special{pa 1000 2200}%
\special{fp}%
\special{pa 1140 2000}%
\special{pa 1000 2140}%
\special{fp}%
\special{pa 1080 2000}%
\special{pa 1000 2080}%
\special{fp}%
%
\special{pn 4}%
\special{pa 860 2400}%
\special{pa 600 2660}%
\special{fp}%
\special{pa 920 2400}%
\special{pa 600 2720}%
\special{fp}%
\special{pa 980 2400}%
\special{pa 600 2780}%
\special{fp}%
\special{pa 1000 2440}%
\special{pa 640 2800}%
\special{fp}%
\special{pa 1000 2500}%
\special{pa 700 2800}%
\special{fp}%
\special{pa 1000 2560}%
\special{pa 760 2800}%
\special{fp}%
\special{pa 1000 2620}%
\special{pa 820 2800}%
\special{fp}%
\special{pa 1000 2680}%
\special{pa 880 2800}%
\special{fp}%
\special{pa 1000 2740}%
\special{pa 940 2800}%
\special{fp}%
\special{pa 800 2400}%
\special{pa 600 2600}%
\special{fp}%
\special{pa 740 2400}%
\special{pa 600 2540}%
\special{fp}%
\special{pa 680 2400}%
\special{pa 600 2480}%
\special{fp}%
%
\special{pn 4}%
\special{pa 6400 1360}%
\special{pa 6160 1600}%
\special{fp}%
\special{pa 6400 1300}%
\special{pa 6100 1600}%
\special{fp}%
\special{pa 6400 1240}%
\special{pa 6040 1600}%
\special{fp}%
\special{pa 6380 1200}%
\special{pa 6000 1580}%
\special{fp}%
\special{pa 6320 1200}%
\special{pa 6000 1520}%
\special{fp}%
\special{pa 6260 1200}%
\special{pa 6000 1460}%
\special{fp}%
\special{pa 6200 1200}%
\special{pa 6000 1400}%
\special{fp}%
\special{pa 6140 1200}%
\special{pa 6000 1340}%
\special{fp}%
\special{pa 6080 1200}%
\special{pa 6000 1280}%
\special{fp}%
\special{pa 6400 1420}%
\special{pa 6220 1600}%
\special{fp}%
\special{pa 6400 1480}%
\special{pa 6280 1600}%
\special{fp}%
\special{pa 6400 1540}%
\special{pa 6340 1600}%
\special{fp}%
%
\special{pn 4}%
\special{pa 6000 2180}%
\special{pa 5780 2400}%
\special{fp}%
\special{pa 6000 2120}%
\special{pa 5720 2400}%
\special{fp}%
\special{pa 6000 2060}%
\special{pa 5660 2400}%
\special{fp}%
\special{pa 5990 2010}%
\special{pa 5610 2390}%
\special{fp}%
\special{pa 5940 2000}%
\special{pa 5600 2340}%
\special{fp}%
\special{pa 5880 2000}%
\special{pa 5600 2280}%
\special{fp}%
\special{pa 5820 2000}%
\special{pa 5600 2220}%
\special{fp}%
\special{pa 5760 2000}%
\special{pa 5600 2160}%
\special{fp}%
\special{pa 5700 2000}%
\special{pa 5600 2100}%
\special{fp}%
\special{pa 5640 2000}%
\special{pa 5600 2040}%
\special{fp}%
\special{pa 6000 2240}%
\special{pa 5840 2400}%
\special{fp}%
\special{pa 6000 2300}%
\special{pa 5900 2400}%
\special{fp}%
\special{pa 6000 2360}%
\special{pa 5960 2400}%
\special{fp}%
%
\special{pn 8}%
\special{sh 0.300}%
\special{pa 1000 3400}%
\special{pa 1200 3400}%
\special{pa 1200 3600}%
\special{pa 1000 3600}%
\special{pa 1000 3400}%
\special{fp}%
%
\special{pn 4}%
\special{pa 1160 3000}%
\special{pa 1000 3160}%
\special{fp}%
\special{pa 1200 3020}%
\special{pa 1020 3200}%
\special{fp}%
\special{pa 1200 3080}%
\special{pa 1080 3200}%
\special{fp}%
\special{pa 1200 3140}%
\special{pa 1140 3200}%
\special{fp}%
\special{pa 1100 3000}%
\special{pa 1000 3100}%
\special{fp}%
\special{pa 1040 3000}%
\special{pa 1000 3040}%
\special{fp}%
\put(14.0000,-30.0000){\makebox(0,0)[lt]{addable box}}%
\put(14.0000,-34.0000){\makebox(0,0)[lt]{removable box}}%
\put(64.0000,-18.0000){\makebox(0,0)[lt]{$0$-addable}}%
%
\special{pn 8}%
\special{pa 6400 1800}%
\special{pa 6200 1400}%
\special{fp}%
\special{sh 1}%
\special{pa 6200 1400}%
\special{pa 6212 1469}%
\special{pa 6224 1448}%
\special{pa 6248 1451}%
\special{pa 6200 1400}%
\special{fp}%
%
\special{pn 8}%
\special{pa 6400 1800}%
\special{pa 5800 2200}%
\special{fp}%
\special{sh 1}%
\special{pa 5800 2200}%
\special{pa 5867 2180}%
\special{pa 5844 2170}%
\special{pa 5844 2146}%
\special{pa 5800 2200}%
\special{fp}%
\end{picture}%
\]
We will use the following notations for the description of the Hayashi realization.
\begin{dfn}
Let $\ell$ be a fixed positive number. Assume that a partition $\lambda$ is obtained by adding a box $x$ to a partition $\mu$. 
\begin{eqnarray*}
N_i^b(\lambda,\mu)&=&\sharp\{y \in \lambda|\text{$y$ is an $i$-addable box below $x$}\} \\
&{}& \quad -\sharp\{y \in \lambda|\text{$y$is an $i$-removable box below $x$}\}, \\
N_i^a(\lambda,\mu)&=&\sharp\{y \in \lambda|\text{$y$ is an $i$-addable box above $x$}\} \\
&{}& \quad -\sharp\{y \in \lambda|\text{$y$ is an $i$-removable box above $x$}\}, \\
N_i(\lambda)&=&\sharp\{y \in \lambda|\text{$y$ is an $i$-addable box}\}-\sharp\{y \in \lambda|\text{$y$ is a $i$-removable box}\}, \\
\end{eqnarray*} 
\end{dfn}
\quad {} \\
Let $|\lambda\rangle$ be a symbol associated to the partition $\lambda \in \mathcal{P}$. The Fock space of $U_q(\widehat{\mf{sl}_\ell})$ is defined
\[
\mathcal{F}=\bigoplus_{\lambda \in \mathcal{P}}\bb{C}(q)|\lambda\rangle.
\]
In \cite{Ha}, T. Hayashi defined the $U_q(\widehat{\mf{sl}_\ell})$-action on the Fock space $\mathcal{F}$.
\begin{thm}[the Hayashi realization of the Fock space] The Fock space $\mathcal{F}$ becomes a $U_q(\widehat{\mf{sl}_n})$-module via the following action;
\begin{eqnarray*}
E_i|\lambda\rangle&=&\sum_{\res_{\ell}(\lambda/\nu)=i}q^{-N_i^a(\nu,\lambda)}|\nu\rangle, \\
F_i|\lambda\rangle&=&\sum_{\res_{\ell}(\mu/\lambda)=i}q^{N_i^b(\lambda,\mu)}|\mu\rangle, \\
K_i|\lambda\rangle&=&q^{N_i(\lambda)}|\lambda\rangle \quad (0 \le i \le \ell-1).
\end{eqnarray*}
Note that the operator $E_i$ removes one box from $\lambda$, and $F_i$ adds one box to $\lambda$. 
\end{thm}
A proof of this theorem can be found in \cite{Ari2}.
\subsubsection{Wedge space}
We will recall the realization of the Fock space as a semi-infinite wedge space. This section is  based on \cite{KMS}.\\
\indent Let $V=\bb{C}^{\ell}$ with basis $v_1, \cdots ,v_\ell$, and $V(z)=V \otimes \bb{C}(q)[z,z^{-1}]$ with basis $u_{j-a\ell}=z^av_j$. Then $U_q(\widehat{\mf{sl}_\ell})$ acts on $V(z)$ by the following way;
\begin{eqnarray*}
E_iu_m&=&\delta(m-1 \equiv i \ \mod{\ell})u_{m-1}, \\
F_iu_m&=&\delta(m \equiv i \ \mod{\ell})u_{m+1}, \\
K_iu_m&=&q^{\delta(m \equiv i \ \mod{\ell})-\delta(m \equiv i+1 \ \mod{\ell})}u_m.
\end{eqnarray*}
This module $V(z)$ is called the evaluation module of $U_q(\widehat{\mf{sl}_\ell})$. 
\quad {} \\
\indent Let $I=( \cdots ,i_2,i_1,i_0)$ be a semi-infinite sequence of integers such that $i_0>i_1>i_2> \cdots$ and $i_k=-k+1$ if $k \gg 0$. Let $u_I$ be semi-infinite wedge product,  
\[
u_I= \cdots \wedge u_{i_2} \wedge u_{i_1} \wedge u_{i_0}.
\]
This wedge product satisfies the following relations; if $k>m$,
\begin{eqnarray*}
u_k \wedge u_m&=&-u_m \wedge u_k \quad (k \equiv m \ \mod\ell), \\
u_k \wedge u_m&=&-qu_m \wedge u_k \\
&{}& \quad +(q^2-1)\left\{
u_{m-i} \wedge u_{k+i}-qu_{m-\ell} \wedge u_{k+\ell}+q^2u_{m-\ell+i} \wedge u_{k+\ell+i}- \cdots 
\right\}\\
&{}& \hspace{80mm} (m-k \equiv i \ \mod{\ell}, 0<i<\ell).
\end{eqnarray*}
Let $\vac_{-k}$ be the $k$-th vacuum vector defined by
\[
\vac_{-k}= \cdots \wedge u_{-(k+2)} \wedge u_{-(k+1)} \wedge u_{-k}.
\]
\begin{dfn}
The Fock space of $U_q(\widehat{\mf{sl}_\ell})$ is defined by
\[
\mathcal{F}=\bigoplus_{I}\bb{C}(q)u_I,
\]
where $I$ runs over the set of semi-infinite increasing sequences of integers such that $i_k=-k+1\ (k \gg 0)$.
\end{dfn}
The Fock space $\mathcal{F}$ becomes a $U_q(\widehat{\mf{sl}_\ell})$-module by the coproduct $\Delta^-$. We have
\begin{eqnarray}
E_i\vac_{-k}&=&0, \label{rv1} \\
F_i\vac_{-k}&=&\left\{
\begin{array}{ll}
\vac_{-k-1} \wedge u_{-k+1} & (i \equiv -k \ \mod\ell) \\
0 & \text{otherwise} 
\end{array}
\right., \label{rv2} \\
K_i\vac_{-k}&=&\left\{
\begin{array}{ll}
q\vac_{-k} & (i \equiv -k \ \mod{\ell}) \\
\vac_{-k} & \text{otherwise}
\end{array}
\right., \label{rv3}
\end{eqnarray}
and the coproduct $\Delta^{-}$.
\begin{prop}
These two realizations are isomorphic as $U_q(\widehat{\mf{sl}_n})$-modules by the one-to-one correspondence
\[
|\lambda=(\lambda_0 \ge \lambda_1 \ge \cdots )\rangle \leftrightarrow \cdots u_{\lambda_2-2} \wedge u_{\lambda_1-1} \wedge u_{\lambda_0}.
\]
\end{prop}
\subsection{Crystal basis : Misra-Miwa's Theorem}
\indent \indent We will use the notion of ``$i$-good box'' for the description of the crystal structure of the Fock space $\mathcal{F}$. 
\begin{dfn}
Let $\ell$ be a fixed positive number, and $\lambda$ be a partition. Reading the $i$-addable boxes and the $i$-removable boxes in $\lambda$ from bottom up, we can obtain the sequence of $A$ and $R$. Next, delete all occurrences of $AR$ from this sequence and keep doing this until no such subsequences remain. Then the $i$-good box in $\lambda$ is the corresponding $i$-removable box to the rightest $R$ in this sequence. 
\end{dfn}
\[
\unitlength 0.1in
\begin{picture}(22.00,25.70)(6.00,-26.00)
%
\special{pn 8}%
\special{pa 800 400}%
\special{pa 1200 400}%
\special{pa 1200 800}%
\special{pa 800 800}%
\special{pa 800 400}%
\special{fp}%
%
\special{pn 8}%
\special{pa 1200 400}%
\special{pa 1600 400}%
\special{pa 1600 800}%
\special{pa 1200 800}%
\special{pa 1200 400}%
\special{fp}%
%
\special{pn 8}%
\special{pa 1600 400}%
\special{pa 2000 400}%
\special{pa 2000 800}%
\special{pa 1600 800}%
\special{pa 1600 400}%
\special{fp}%
%
\special{pn 8}%
\special{pa 2000 400}%
\special{pa 2400 400}%
\special{pa 2400 800}%
\special{pa 2000 800}%
\special{pa 2000 400}%
\special{fp}%
%
\special{pn 8}%
\special{pa 800 800}%
\special{pa 1200 800}%
\special{pa 1200 1200}%
\special{pa 800 1200}%
\special{pa 800 800}%
\special{fp}%
%
\special{pn 8}%
\special{pa 1200 800}%
\special{pa 1600 800}%
\special{pa 1600 1200}%
\special{pa 1200 1200}%
\special{pa 1200 800}%
\special{fp}%
%
\special{pn 8}%
\special{pa 1600 800}%
\special{pa 2000 800}%
\special{pa 2000 1200}%
\special{pa 1600 1200}%
\special{pa 1600 800}%
\special{fp}%
%
\special{pn 8}%
\special{pa 800 1200}%
\special{pa 1200 1200}%
\special{pa 1200 1600}%
\special{pa 800 1600}%
\special{pa 800 1200}%
\special{fp}%
\put(10.0000,-6.0000){\makebox(0,0){$0$}}%
\put(14.0000,-10.0000){\makebox(0,0){$0$}}%
\put(14.0000,-6.0000){\makebox(0,0){$1$}}%
\put(18.0000,-6.0000){\makebox(0,0){$2$}}%
\put(21.9000,-6.0000){\makebox(0,0){$0$}}%
\put(18.0000,-10.0000){\makebox(0,0){$1$}}%
\put(10.0000,-10.0000){\makebox(0,0){$2$}}%
\put(10.0000,-14.0000){\makebox(0,0){$1$}}%
%
\special{pn 8}%
\special{pa 1200 1200}%
\special{pa 1600 1200}%
\special{pa 1600 1600}%
\special{pa 1200 1600}%
\special{pa 1200 1200}%
\special{fp}%
\put(14.0000,-14.0000){\makebox(0,0){$2$}}%
%
\special{pn 8}%
\special{pa 1200 1600}%
\special{pa 1600 1600}%
\special{pa 1600 2000}%
\special{pa 1200 2000}%
\special{pa 1200 1600}%
\special{fp}%
\put(14.0000,-18.0000){\makebox(0,0){$1$}}%
%
\special{pn 8}%
\special{pa 800 1600}%
\special{pa 1200 1600}%
\special{pa 1200 2000}%
\special{pa 800 2000}%
\special{pa 800 1600}%
\special{fp}%
\put(10.0000,-18.0000){\makebox(0,0){$0$}}%
%
\special{pn 8}%
\special{pa 800 2000}%
\special{pa 1200 2000}%
\special{pa 1200 2400}%
\special{pa 800 2400}%
\special{pa 800 2000}%
\special{fp}%
\put(10.0000,-22.0000){\makebox(0,0){$2$}}%
\put(6.0000,-2.0000){\makebox(0,0)[lb]{$\lambda=(4,3,2,2,1), \ \ell=3$}}%
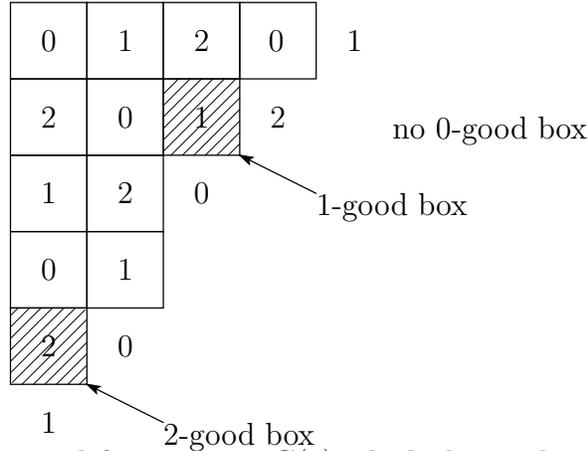
\put(14.0000,-22.0000){\makebox(0,0){$0$}}%
\put(18.0000,-14.0000){\makebox(0,0){$0$}}%
\put(10.0000,-26.0000){\makebox(0,0){$1$}}%
\put(22.0000,-10.0000){\makebox(0,0){$2$}}%
\put(26.0000,-6.0000){\makebox(0,0){$1$}}%
%
\special{pn 4}%
\special{pa 2000 960}%
\special{pa 1760 1200}%
\special{fp}%
\special{pa 2000 900}%
\special{pa 1700 1200}%
\special{fp}%
\special{pa 2000 840}%
\special{pa 1640 1200}%
\special{fp}%
\special{pa 1980 800}%
\special{pa 1600 1180}%
\special{fp}%
\special{pa 1920 800}%
\special{pa 1600 1120}%
\special{fp}%
\special{pa 1860 800}%
\special{pa 1600 1060}%
\special{fp}%
\special{pa 1800 800}%
\special{pa 1600 1000}%
\special{fp}%
\special{pa 1740 800}%
\special{pa 1600 940}%
\special{fp}%
\special{pa 1680 800}%
\special{pa 1600 880}%
\special{fp}%
\special{pa 2000 1020}%
\special{pa 1820 1200}%
\special{fp}%
\special{pa 2000 1080}%
\special{pa 1880 1200}%
\special{fp}%
\special{pa 2000 1140}%
\special{pa 1940 1200}%
\special{fp}%
%
\special{pn 4}%
\special{pa 1200 2240}%
\special{pa 1040 2400}%
\special{fp}%
\special{pa 1200 2180}%
\special{pa 980 2400}%
\special{fp}%
\special{pa 1200 2120}%
\special{pa 920 2400}%
\special{fp}%
\special{pa 1200 2060}%
\special{pa 860 2400}%
\special{fp}%
\special{pa 1190 2010}%
\special{pa 810 2390}%
\special{fp}%
\special{pa 1140 2000}%
\special{pa 800 2340}%
\special{fp}%
\special{pa 1080 2000}%
\special{pa 800 2280}%
\special{fp}%
\special{pa 1020 2000}%
\special{pa 800 2220}%
\special{fp}%
\special{pa 960 2000}%
\special{pa 800 2160}%
\special{fp}%
\special{pa 900 2000}%
\special{pa 800 2100}%
\special{fp}%
\special{pa 840 2000}%
\special{pa 800 2040}%
\special{fp}%
\special{pa 1200 2300}%
\special{pa 1100 2400}%
\special{fp}%
\special{pa 1200 2360}%
\special{pa 1160 2400}%
\special{fp}%
\put(28.0000,-10.0000){\makebox(0,0)[lt]{no $0$-good box}}%
%
\special{pn 8}%
\special{pa 2400 1400}%
\special{pa 2000 1200}%
\special{fp}%
\special{sh 1}%
\special{pa 2000 1200}%
\special{pa 2051 1248}%
\special{pa 2048 1224}%
\special{pa 2069 1212}%
\special{pa 2000 1200}%
\special{fp}%
\put(24.0000,-14.0000){\makebox(0,0)[lt]{$1$-good box}}%
\put(16.0000,-26.0000){\makebox(0,0)[lt]{$2$-good box}}%
%
\special{pn 8}%
\special{pa 1600 2600}%
\special{pa 1200 2400}%
\special{fp}%
\special{sh 1}%
\special{pa 1200 2400}%
\special{pa 1251 2448}%
\special{pa 1248 2424}%
\special{pa 1269 2412}%
\special{pa 1200 2400}%
\special{fp}%
\end{picture}%
\]
\indent Let $R$ be the subring of rational functions in $\bb{C}(q)$ which do not have a pole at $0$. Let
\[
L=\bigoplus_{\lambda \in \mathcal{P}}R|\lambda\rangle, \quad  B=\{|\lambda \rangle \ (\mod{qL})\}.
\] 
\begin{thm}[Misra-Miwa \cite{MM}]\label{MM} The $(L,B)$ is a crystal basis of $\mathcal{F}$ by the following action of Kashiwara operators $\widetilde{e_i}, \widetilde{f_i} \ (1 \le i \le \ell-1)$; \\
\indent (1) \ If a partition $\lambda$ has no $i$-good box, then $\widetilde{e_i}|\lambda \rangle=0 \ (\mod{qL})$. \\
\indent (2) \ If $x$ is a $i$-good box of $\lambda$ and $\mu=\lambda \backslash \{x\}$,
\[
\widetilde{e_i}|\lambda\rangle=|\mu\rangle \ (\mod{qL}), \quad \widetilde{f_i}|\mu\rangle=|\lambda\rangle. \ (\mod{qL}).
\]
\indent (3) \ If a partition $\mu$ has no $i$-addable box $x$ which is $i$-good box in $\mu \cup \{x\}$, then $\widetilde{f_i}|\mu\rangle=0 \ (\mod{qL})$. 
\end{thm}
A proof of this theorem can be found in \cite{Ari2}.
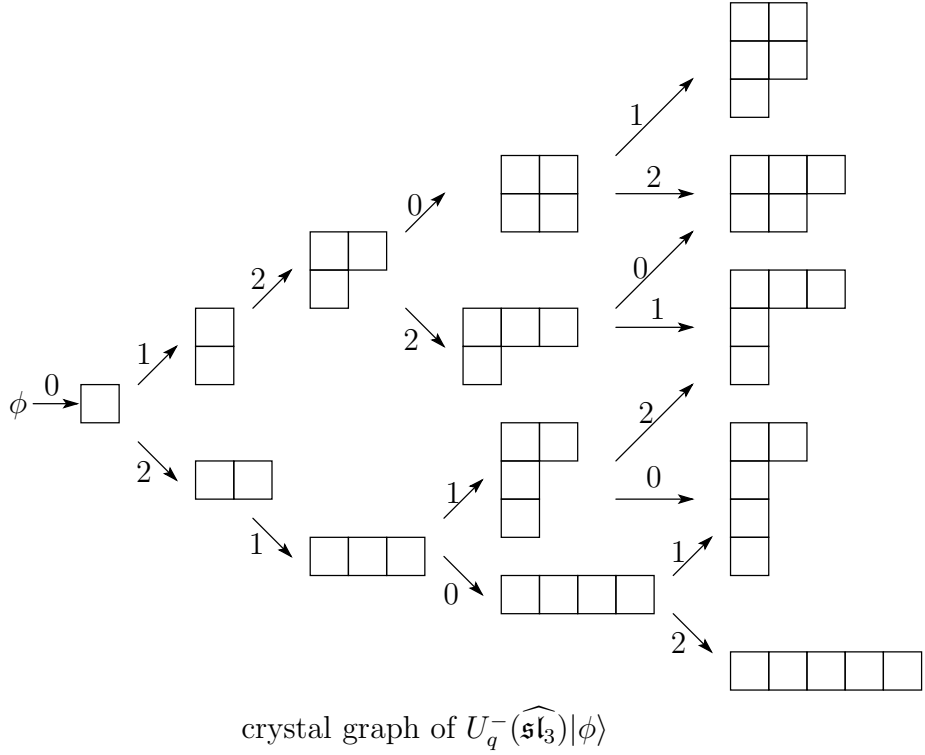
\[
\unitlength 0.1in
\begin{picture}(52.55,37.15)(3.45,-39.15)
\put(9.2000,-23.8000){\makebox(0,0)[rb]{$\phi$}}%
%
\special{pn 8}%
\special{pa 1200 2200}%
\special{pa 1400 2200}%
\special{pa 1400 2400}%
\special{pa 1200 2400}%
\special{pa 1200 2200}%
\special{fp}%
%
\special{pn 8}%
\special{pa 1800 2000}%
\special{pa 2000 2000}%
\special{pa 2000 2200}%
\special{pa 1800 2200}%
\special{pa 1800 2000}%
\special{fp}%
%
\special{pn 8}%
\special{pa 1800 1800}%
\special{pa 2000 1800}%
\special{pa 2000 2000}%
\special{pa 1800 2000}%
\special{pa 1800 1800}%
\special{fp}%
%
\special{pn 8}%
\special{pa 1800 2600}%
\special{pa 2000 2600}%
\special{pa 2000 2800}%
\special{pa 1800 2800}%
\special{pa 1800 2600}%
\special{fp}%
%
\special{pn 8}%
\special{pa 2000 2600}%
\special{pa 2200 2600}%
\special{pa 2200 2800}%
\special{pa 2000 2800}%
\special{pa 2000 2600}%
\special{fp}%
%
\special{pn 8}%
\special{pa 2400 1400}%
\special{pa 2600 1400}%
\special{pa 2600 1600}%
\special{pa 2400 1600}%
\special{pa 2400 1400}%
\special{fp}%
%
\special{pn 8}%
\special{pa 2600 1400}%
\special{pa 2800 1400}%
\special{pa 2800 1600}%
\special{pa 2600 1600}%
\special{pa 2600 1400}%
\special{fp}%
%
\special{pn 8}%
\special{pa 2400 1600}%
\special{pa 2600 1600}%
\special{pa 2600 1800}%
\special{pa 2400 1800}%
\special{pa 2400 1600}%
\special{fp}%
%
\special{pn 8}%
\special{pa 2400 3000}%
\special{pa 2600 3000}%
\special{pa 2600 3200}%
\special{pa 2400 3200}%
\special{pa 2400 3000}%
\special{fp}%
%
\special{pn 8}%
\special{pa 2600 3000}%
\special{pa 2800 3000}%
\special{pa 2800 3200}%
\special{pa 2600 3200}%
\special{pa 2600 3000}%
\special{fp}%
%
\special{pn 8}%
\special{pa 2800 3000}%
\special{pa 3000 3000}%
\special{pa 3000 3200}%
\special{pa 2800 3200}%
\special{pa 2800 3000}%
\special{fp}%
%
\special{pn 8}%
\special{pa 3400 1000}%
\special{pa 3600 1000}%
\special{pa 3600 1200}%
\special{pa 3400 1200}%
\special{pa 3400 1000}%
\special{fp}%
%
\special{pn 8}%
\special{pa 3600 1000}%
\special{pa 3800 1000}%
\special{pa 3800 1200}%
\special{pa 3600 1200}%
\special{pa 3600 1000}%
\special{fp}%
%
\special{pn 8}%
\special{pa 3400 1200}%
\special{pa 3600 1200}%
\special{pa 3600 1400}%
\special{pa 3400 1400}%
\special{pa 3400 1200}%
\special{fp}%
%
\special{pn 8}%
\special{pa 3600 1200}%
\special{pa 3800 1200}%
\special{pa 3800 1400}%
\special{pa 3600 1400}%
\special{pa 3600 1200}%
\special{fp}%
%
\special{pn 8}%
\special{pa 3200 1800}%
\special{pa 3400 1800}%
\special{pa 3400 2000}%
\special{pa 3200 2000}%
\special{pa 3200 1800}%
\special{fp}%
%
\special{pn 8}%
\special{pa 3400 1800}%
\special{pa 3600 1800}%
\special{pa 3600 2000}%
\special{pa 3400 2000}%
\special{pa 3400 1800}%
\special{fp}%
%
\special{pn 8}%
\special{pa 3200 2000}%
\special{pa 3400 2000}%
\special{pa 3400 2200}%
\special{pa 3200 2200}%
\special{pa 3200 2000}%
\special{fp}%
%
\special{pn 8}%
\special{pa 3600 1800}%
\special{pa 3800 1800}%
\special{pa 3800 2000}%
\special{pa 3600 2000}%
\special{pa 3600 1800}%
\special{fp}%
%
\special{pn 8}%
\special{pa 3400 2400}%
\special{pa 3600 2400}%
\special{pa 3600 2600}%
\special{pa 3400 2600}%
\special{pa 3400 2400}%
\special{fp}%
%
\special{pn 8}%
\special{pa 3600 2400}%
\special{pa 3800 2400}%
\special{pa 3800 2600}%
\special{pa 3600 2600}%
\special{pa 3600 2400}%
\special{fp}%
%
\special{pn 8}%
\special{pa 3400 2600}%
\special{pa 3600 2600}%
\special{pa 3600 2800}%
\special{pa 3400 2800}%
\special{pa 3400 2600}%
\special{fp}%
%
\special{pn 8}%
\special{pa 3400 2800}%
\special{pa 3600 2800}%
\special{pa 3600 3000}%
\special{pa 3400 3000}%
\special{pa 3400 2800}%
\special{fp}%
%
\special{pn 8}%
\special{pa 3400 3200}%
\special{pa 3600 3200}%
\special{pa 3600 3400}%
\special{pa 3400 3400}%
\special{pa 3400 3200}%
\special{fp}%
%
\special{pn 8}%
\special{pa 3600 3200}%
\special{pa 3800 3200}%
\special{pa 3800 3400}%
\special{pa 3600 3400}%
\special{pa 3600 3200}%
\special{fp}%
%
\special{pn 8}%
\special{pa 3800 3200}%
\special{pa 4000 3200}%
\special{pa 4000 3400}%
\special{pa 3800 3400}%
\special{pa 3800 3200}%
\special{fp}%
%
\special{pn 8}%
\special{pa 4000 3200}%
\special{pa 4200 3200}%
\special{pa 4200 3400}%
\special{pa 4000 3400}%
\special{pa 4000 3200}%
\special{fp}%
%
\special{pn 8}%
\special{pa 4600 2400}%
\special{pa 4800 2400}%
\special{pa 4800 2600}%
\special{pa 4600 2600}%
\special{pa 4600 2400}%
\special{fp}%
%
\special{pn 8}%
\special{pa 4800 2400}%
\special{pa 5000 2400}%
\special{pa 5000 2600}%
\special{pa 4800 2600}%
\special{pa 4800 2400}%
\special{fp}%
%
\special{pn 8}%
\special{pa 4600 2600}%
\special{pa 4800 2600}%
\special{pa 4800 2800}%
\special{pa 4600 2800}%
\special{pa 4600 2600}%
\special{fp}%
%
\special{pn 8}%
\special{pa 4600 2800}%
\special{pa 4800 2800}%
\special{pa 4800 3000}%
\special{pa 4600 3000}%
\special{pa 4600 2800}%
\special{fp}%
%
\special{pn 8}%
\special{pa 4600 3000}%
\special{pa 4800 3000}%
\special{pa 4800 3200}%
\special{pa 4600 3200}%
\special{pa 4600 3000}%
\special{fp}%
%
\special{pn 8}%
\special{pa 4600 3600}%
\special{pa 4800 3600}%
\special{pa 4800 3800}%
\special{pa 4600 3800}%
\special{pa 4600 3600}%
\special{fp}%
%
\special{pn 8}%
\special{pa 4800 3600}%
\special{pa 5000 3600}%
\special{pa 5000 3800}%
\special{pa 4800 3800}%
\special{pa 4800 3600}%
\special{fp}%
%
\special{pn 8}%
\special{pa 5000 3600}%
\special{pa 5200 3600}%
\special{pa 5200 3800}%
\special{pa 5000 3800}%
\special{pa 5000 3600}%
\special{fp}%
%
\special{pn 8}%
\special{pa 5200 3600}%
\special{pa 5400 3600}%
\special{pa 5400 3800}%
\special{pa 5200 3800}%
\special{pa 5200 3600}%
\special{fp}%
%
\special{pn 8}%
\special{pa 5400 3600}%
\special{pa 5600 3600}%
\special{pa 5600 3800}%
\special{pa 5400 3800}%
\special{pa 5400 3600}%
\special{fp}%
%
\special{pn 8}%
\special{pa 4600 1600}%
\special{pa 4800 1600}%
\special{pa 4800 1800}%
\special{pa 4600 1800}%
\special{pa 4600 1600}%
\special{fp}%
%
\special{pn 8}%
\special{pa 4800 1600}%
\special{pa 5000 1600}%
\special{pa 5000 1800}%
\special{pa 4800 1800}%
\special{pa 4800 1600}%
\special{fp}%
%
\special{pn 8}%
\special{pa 4600 1800}%
\special{pa 4800 1800}%
\special{pa 4800 2000}%
\special{pa 4600 2000}%
\special{pa 4600 1800}%
\special{fp}%
%
\special{pn 8}%
\special{pa 4600 2000}%
\special{pa 4800 2000}%
\special{pa 4800 2200}%
\special{pa 4600 2200}%
\special{pa 4600 2000}%
\special{fp}%
%
\special{pn 8}%
\special{pa 5000 1600}%
\special{pa 5200 1600}%
\special{pa 5200 1800}%
\special{pa 5000 1800}%
\special{pa 5000 1600}%
\special{fp}%
%
\special{pn 8}%
\special{pa 4600 1000}%
\special{pa 4800 1000}%
\special{pa 4800 1200}%
\special{pa 4600 1200}%
\special{pa 4600 1000}%
\special{fp}%
%
\special{pn 8}%
\special{pa 4800 1000}%
\special{pa 5000 1000}%
\special{pa 5000 1200}%
\special{pa 4800 1200}%
\special{pa 4800 1000}%
\special{fp}%
%
\special{pn 8}%
\special{pa 4600 1200}%
\special{pa 4800 1200}%
\special{pa 4800 1400}%
\special{pa 4600 1400}%
\special{pa 4600 1200}%
\special{fp}%
%
\special{pn 8}%
\special{pa 4800 1200}%
\special{pa 5000 1200}%
\special{pa 5000 1400}%
\special{pa 4800 1400}%
\special{pa 4800 1200}%
\special{fp}%
%
\special{pn 8}%
\special{pa 5000 1000}%
\special{pa 5200 1000}%
\special{pa 5200 1200}%
\special{pa 5000 1200}%
\special{pa 5000 1000}%
\special{fp}%
%
\special{pn 8}%
\special{pa 4600 200}%
\special{pa 4800 200}%
\special{pa 4800 400}%
\special{pa 4600 400}%
\special{pa 4600 200}%
\special{fp}%
%
\special{pn 8}%
\special{pa 4800 200}%
\special{pa 5000 200}%
\special{pa 5000 400}%
\special{pa 4800 400}%
\special{pa 4800 200}%
\special{fp}%
%
\special{pn 8}%
\special{pa 4600 400}%
\special{pa 4800 400}%
\special{pa 4800 600}%
\special{pa 4600 600}%
\special{pa 4600 400}%
\special{fp}%
%
\special{pn 8}%
\special{pa 4800 400}%
\special{pa 5000 400}%
\special{pa 5000 600}%
\special{pa 4800 600}%
\special{pa 4800 400}%
\special{fp}%
%
\special{pn 8}%
\special{pa 4600 600}%
\special{pa 4800 600}%
\special{pa 4800 800}%
\special{pa 4600 800}%
\special{pa 4600 600}%
\special{fp}%
%
\special{pn 8}%
\special{pa 950 2300}%
\special{pa 1180 2300}%
\special{fp}%
\special{sh 1}%
\special{pa 1180 2300}%
\special{pa 1113 2280}%
\special{pa 1127 2300}%
\special{pa 1113 2320}%
\special{pa 1180 2300}%
\special{fp}%
%
\special{pn 8}%
\special{pa 4300 3200}%
\special{pa 4500 3000}%
\special{fp}%
\special{sh 1}%
\special{pa 4500 3000}%
\special{pa 4439 3033}%
\special{pa 4462 3038}%
\special{pa 4467 3061}%
\special{pa 4500 3000}%
\special{fp}%
%
\special{pn 8}%
\special{pa 1500 2200}%
\special{pa 1700 2000}%
\special{fp}%
\special{sh 1}%
\special{pa 1700 2000}%
\special{pa 1639 2033}%
\special{pa 1662 2038}%
\special{pa 1667 2061}%
\special{pa 1700 2000}%
\special{fp}%
%
\special{pn 8}%
\special{pa 2100 1800}%
\special{pa 2300 1600}%
\special{fp}%
\special{sh 1}%
\special{pa 2300 1600}%
\special{pa 2239 1633}%
\special{pa 2262 1638}%
\special{pa 2267 1661}%
\special{pa 2300 1600}%
\special{fp}%
%
\special{pn 8}%
\special{pa 2900 1400}%
\special{pa 3100 1200}%
\special{fp}%
\special{sh 1}%
\special{pa 3100 1200}%
\special{pa 3039 1233}%
\special{pa 3062 1238}%
\special{pa 3067 1261}%
\special{pa 3100 1200}%
\special{fp}%
%
\special{pn 8}%
\special{pa 4000 1000}%
\special{pa 4400 600}%
\special{fp}%
\special{sh 1}%
\special{pa 4400 600}%
\special{pa 4339 633}%
\special{pa 4362 638}%
\special{pa 4367 661}%
\special{pa 4400 600}%
\special{fp}%
%
\special{pn 8}%
\special{pa 4000 1200}%
\special{pa 4400 1200}%
\special{fp}%
\special{sh 1}%
\special{pa 4400 1200}%
\special{pa 4333 1180}%
\special{pa 4347 1200}%
\special{pa 4333 1220}%
\special{pa 4400 1200}%
\special{fp}%
%
\special{pn 8}%
\special{pa 4000 1900}%
\special{pa 4400 1900}%
\special{fp}%
\special{sh 1}%
\special{pa 4400 1900}%
\special{pa 4333 1880}%
\special{pa 4347 1900}%
\special{pa 4333 1920}%
\special{pa 4400 1900}%
\special{fp}%
%
\special{pn 8}%
\special{pa 4000 2800}%
\special{pa 4400 2800}%
\special{fp}%
\special{sh 1}%
\special{pa 4400 2800}%
\special{pa 4333 2780}%
\special{pa 4347 2800}%
\special{pa 4333 2820}%
\special{pa 4400 2800}%
\special{fp}%
%
\special{pn 8}%
\special{pa 4000 2600}%
\special{pa 4400 2200}%
\special{fp}%
\special{sh 1}%
\special{pa 4400 2200}%
\special{pa 4339 2233}%
\special{pa 4362 2238}%
\special{pa 4367 2261}%
\special{pa 4400 2200}%
\special{fp}%
%
\special{pn 8}%
\special{pa 4000 1800}%
\special{pa 4400 1400}%
\special{fp}%
\special{sh 1}%
\special{pa 4400 1400}%
\special{pa 4339 1433}%
\special{pa 4362 1438}%
\special{pa 4367 1461}%
\special{pa 4400 1400}%
\special{fp}%
%
\special{pn 8}%
\special{pa 1500 2500}%
\special{pa 1700 2700}%
\special{fp}%
\special{sh 1}%
\special{pa 1700 2700}%
\special{pa 1667 2639}%
\special{pa 1662 2662}%
\special{pa 1639 2667}%
\special{pa 1700 2700}%
\special{fp}%
%
\special{pn 8}%
\special{pa 2100 2900}%
\special{pa 2300 3100}%
\special{fp}%
\special{sh 1}%
\special{pa 2300 3100}%
\special{pa 2267 3039}%
\special{pa 2262 3062}%
\special{pa 2239 3067}%
\special{pa 2300 3100}%
\special{fp}%
%
\special{pn 8}%
\special{pa 3100 3100}%
\special{pa 3300 3300}%
\special{fp}%
\special{sh 1}%
\special{pa 3300 3300}%
\special{pa 3267 3239}%
\special{pa 3262 3262}%
\special{pa 3239 3267}%
\special{pa 3300 3300}%
\special{fp}%
%
\special{pn 8}%
\special{pa 3100 2900}%
\special{pa 3300 2700}%
\special{fp}%
\special{sh 1}%
\special{pa 3300 2700}%
\special{pa 3239 2733}%
\special{pa 3262 2738}%
\special{pa 3267 2761}%
\special{pa 3300 2700}%
\special{fp}%
%
\special{pn 8}%
\special{pa 2900 1800}%
\special{pa 3100 2000}%
\special{fp}%
\special{sh 1}%
\special{pa 3100 2000}%
\special{pa 3067 1939}%
\special{pa 3062 1962}%
\special{pa 3039 1967}%
\special{pa 3100 2000}%
\special{fp}%
%
\special{pn 8}%
\special{pa 4300 3400}%
\special{pa 4500 3600}%
\special{fp}%
\special{sh 1}%
\special{pa 4500 3600}%
\special{pa 4467 3539}%
\special{pa 4462 3562}%
\special{pa 4439 3567}%
\special{pa 4500 3600}%
\special{fp}%
\put(10.5000,-22.0000){\makebox(0,0){$0$}}%
\put(15.7000,-20.9000){\makebox(0,0)[rb]{$1$}}%
\put(21.7000,-17.0000){\makebox(0,0)[rb]{$2$}}%
\put(29.9000,-13.1000){\makebox(0,0)[rb]{$0$}}%
\put(41.5000,-8.4000){\makebox(0,0)[rb]{$1$}}%
\put(15.7000,-26.0000){\makebox(0,0)[rt]{$2$}}%
\put(21.6000,-29.8000){\makebox(0,0)[rt]{$1$}}%
\put(31.8000,-32.4000){\makebox(0,0)[rt]{$0$}}%
\put(43.7000,-35.0000){\makebox(0,0)[rt]{$2$}}%
\put(42.0000,-11.1000){\makebox(0,0){$2$}}%
\put(42.2000,-18.0000){\makebox(0,0){$1$}}%
\put(42.0000,-26.8000){\makebox(0,0){$0$}}%
\put(31.9000,-28.2000){\makebox(0,0)[rb]{$1$}}%
\put(43.7000,-31.3000){\makebox(0,0)[rb]{$1$}}%
\put(42.0000,-24.0000){\makebox(0,0)[rb]{$2$}}%
\put(41.7000,-16.4000){\makebox(0,0)[rb]{$0$}}%
\put(29.7000,-19.1000){\makebox(0,0)[rt]{$2$}}%
\put(30.0000,-40.0000){\makebox(0,0){crystal graph of $U_q^{-}(\widehat{\mf{sl}_3})|\phi\rangle$}}%
\end{picture}%
\]
\subsection{Global basis}
In this section, we will introduce the lower and upper global basis of $\mathcal{F}$. We consider the Fock space $\mathcal{F}$ as a wedge space.\\
\indent In \cite{KMS}, the operator $B_k \ (k \in \bb{Z}, k \neq 0)$ on $\mathcal{F}$ is defined by the following;
\[
B_ku_I=(\cdots \wedge u_{i_2} \wedge u_{i_1} \wedge u_{i_0-\ell{k}})+(\cdots \wedge u_{i_2} \wedge u_{i_1-\ell{k}} \wedge u_{i_0})+(\cdots \wedge u_{i_2-\ell{k}} \wedge u_{i_1} \wedge u_{i_0})+\cdots.
\]
Note that $ \cdots u_{i_\nu-\ell{k}} \wedge \cdots \wedge u_{i_0}=0$ for $\nu \gg 0$.
\subsubsection{Lower global basis}
First, we will introduce the bar involution on $\mathcal{F}$. 
\begin{prop}
There exists a unique bar involution ${}^-:\mathcal{F} \to \mathcal{F}$ satisfying the following three properties; \\
\indent (1) \ $\overline{F_iv}=F_i\overline{v} \quad (v \in \mathcal{F}, \ 0 \le i \le \ell-1)$, \\
\indent (2) \ $\overline{B_kv}=B_k\overline{v} \quad (k>0)$, \\
\indent (3) \ $\overline{\vac_0}=\vac_0$. \\
\indent (4) \ $\overline{qv}=q^{-1}\overline{v}$.
\end{prop}
\begin{thm}\label{thm510}
There exists a unique basis 
\[
\{G^{low}(\mu) \in \mathcal{F}|\mu \in \mathcal{P}\}
\]
on $\mathcal{F}$ satisfying the following two properties; 
\begin{enumerate}[(1)]
\item (``bar invariance'') 
\[
\overline{G^{low}(\mu)}=G^{low}(\mu).
\]
\item If $\mu$ is a partiton of $n$, there exists some polynomials $d_{\lambda\mu}(q) \in q\bb{Z}[q]$, then
\[
G^{low}(\mu)=|\mu\rangle+\sum_{\mu \triangleleft \lambda \in \mathcal{P}_n}d_{\lambda\mu}(q)|\lambda\rangle,
\]
where the ordering $\triangleright$ is the dominance ordering of partitions. 
\end{enumerate}
\indent This $\{G^{low}(\mu)\}$ is called the lower global basis of $\mathcal{F}$.
\end{thm}
\begin{exa}\label{exalow} We will calculate the lower global basis for $\ell=2$.\\
(0) \  $\vac_0=|\phi\rangle$ is bar-invariant. Thus
\[
G(\phi)=|\phi\rangle.
\]
(1) \ Since $F_0$ is bar-invariant, 
\[
F_0\vac_0=vac_{-1} \wedge u_1=|(1)\rangle
\]
is bar-invariant. Thus 
\[
G^{low}((1))=|(1)\rangle.
\]
(2) \ Since $F_1$ is bar-invariant, 
\begin{eqnarray*}
F_1(\vac_{-1} \wedge u_1)&=&(F_1\vac_{-1}) \wedge u_1+(K_1\vac_{-1}) \wedge (F_1u_1) \\
&=&\vac_{-2} \wedge u_0 \wedge u_1+q\vac_{-1} \wedge u_2 \\
&=&|(1^2)\rangle+q|(2)\rangle
\end{eqnarray*}
is bar-invariant. Thus
\[
G^{low}((1^2))=|(1^2)\rangle+q|(2)\rangle.
\]
It is clear that $G^{low}((2))=|(2)\rangle$. Therefore 
\[
D=(d_{\lambda\mu})=\left(
\begin{array}{cc}
1 & 0 \\
q & 1 
\end{array}
\right).
\]
(3) \ Since $F_0$ is bar-invariant,
\begin{eqnarray*}
&{}&F_0(\vac_{-2} \wedge u_0 \wedge u_1+q\vac_{-1} \wedge u_2) \\
&=&(F_0vac_{-2}) \wedge u_0 \wedge u_1+(K_0\vac_{-2}) \wedge F_0(u_1 \wedge u_0)+q(F_0\vac_{-1}) \wedge u_2+q(K_0\vac_{-1}) \wedge F_0u_2 \\
&=&vac_{-3} \wedge u_{-1} \wedge u_{0} \wedge u_1+q\vac_{-2} \wedge F_0(u_0 \wedge u_1)+q\vac_{-1} \wedge u_3 \\
&=&vac_{-3} \wedge u_{-1} \wedge u_{0} \wedge u_1+q\vac_{-1} \wedge u_3 \\
&=&|(1^3)\rangle+q|(3)\rangle
\end{eqnarray*}
is bar-invariant. Thus
\[
G^{low}((1^3))=|(1^3)\rangle+q|(3)\rangle.
\]
Since $F_1$ is bar-invariant,
\begin{eqnarray*}
F_1(\vac_{-2} \wedge u_0 \wedge u_1+q\vac_{-1} \wedge u_2)=(q+q^{-1})(\vac_{-2} \wedge u_0 \wedge u_2)
\end{eqnarray*}
are bar-invariant. Moreover, $\vac_{-2} \wedge u_0 \wedge u_2=|(2,1)\rangle$ is bar-invariant because $(q+q^{-1}$ is bar-invariant. Thus
\[
G^{low}((2,1))=|(2,1)\rangle.
\]
It is clear that $G^{low}((3))=|(3)\rangle$. Therefore
\[
D=(d_{\lambda\mu})=\left(
\begin{array}{ccc}
1 & 0 & 0 \\
0 & 1 & 0 \\
q & 0 & 1
\end{array}
\right).
\]
\end{exa}

\subsubsection{Upper global basis}
Let $\{\langle \lambda|\}_{\lambda \in \mathcal{P}}$ be the dual basis of $\{|\lambda\rangle\}_{\lambda \in \mathcal{P}}$ with respect to the scalar product $\langle \lambda|\mu\rangle=\delta_{\lambda\mu}$. The $U_q(\widehat{\mf{sl}_\ell})$-module
\[
\mathcal{F}^\vee=\bigoplus_{\lambda \in \mathcal{P}}\bb{C}(q)\langle\lambda|
\] 
is isomorphic to the wedge space 
\[
\bigoplus_{I} \bb{C}(q) u_I
\]
with the action of $U_q(\widehat{\mf{sl}_n})$ defined by (\ref{rv1}),(\ref{rv2}),(\ref{rv3}) and coproduct $\Delta^+$. 
\begin{prop}
There exists a unique bar involution ${}^-:\mathcal{F}^\vee \to \mathcal{F}^\vee$ satisfying the following three properties; \\
\indent (1) \ $\overline{F_iv}=F_i\overline{v} \quad (v \in \mathcal{F}^\vee, \ 0 \le i \le \ell-1)$, \\
\indent (2) \ $\overline{B_kv}=B_k\overline{v} \quad (k<0)$, \\
\indent (3) \ $\overline{\vac_0}=\vac_0$. \\
\indent (4) \ $\overline{qv}=q^{-1}\overline{v}$
\end{prop}
\begin{thm}
There exists a unique basis 
\[
\{G^{up}(\mu) \in \mathcal{F}^\vee|\mu \in \mathcal{P}\}
\]
on $\mathcal{F}^\vee$ satisfying the following two properties; 
\begin{enumerate}[(1)]
\item (``bar invariance'')
\[
\overline{G^{up}(\mu)}=G^{up}(\mu).
\]
\item If $\lambda$ is a partiton of $n$, there exists some polynomials $d'_{\lambda\mu}(q) \in q\bb{Z}[q]$ then
\[
\langle\lambda|=G^{up}(\lambda)+\sum_{\lambda \triangleright \mu \in \mathcal{P}_n}d_{\lambda\mu}(q)G^{up}(\mu).
\]
Moreover $d'_{\lambda\mu}(q)$ is equal to $d_{\lambda\mu}(q)$ in Theorem \ref{thm510}. Especially the basis $\{G^{up}(\mu)\}$ is the dual basis of $\{G^{low}(\mu)\}$. 
\end{enumerate}
\indent This $\{G^{up}(\mu)\}$ is called the upper global basis of $\mathcal{F}^\vee$.
\end{thm}
\begin{exa}\label{exaup} We will calculate the upper global basis for $\ell=2$.\\
(0) \ Since $\vac_0$ is bar-invariant, $G^{up}(\phi)=\langle \phi|$. \\
(1) \ Since $F_0(\vac_0)=\vac_{-1} \wedge u_1$ is bar-invariant, $G^{up}((1))=\langle(1)|$. \\
(2) \ Note that the action of $F_1$ is defined by the coproduct $\Delta^+$. We obtain
\begin{eqnarray}
F_1(\vac_{-1} \wedge u_1)&=&F_1\vac_{-1} \wedge K_1^{-1}u_1+\vac_{-1} \wedge F_1u_1 \nonumber \\
&=&q^{-1}\vac_{-2} \wedge u_0 \wedge u_1+\vac_{-1} \wedge u_2. \label{up21}
\end{eqnarray}
This is bar-invariant, but not $G^{up}(\vac_{-1} \wedge u_2)=G^{up}((2))$ because the right hand side does not satisfy the condition (2) of the upper global basis, i.e. $q^{-1} \notin q\bb{Z}[q]$. \\
\indent $B_{-1}$ is bar-invariant. Thus 
\begin{eqnarray}
B_{-1}\vac_0&=&\vac_{-2} \wedge u_1 \wedge u_0+\vac_{-2} \wedge u_{-1} \wedge u_2 \nonumber \\
&=&-q\vac_{-2} \wedge u_0 \wedge u_1+\vac_{-1} \wedge u_2 \label{up22}
\end{eqnarray}
is bar-invariant. 
\begin{eqnarray}
(\ref{up21})-(\ref{up22})=(q+q^{-1})\vac_{-2} \wedge u_0 \wedge u_1 \label{up23}
\end{eqnarray}
is bar-invariant.
\[
(\ref{up21})-(\ref{up23})=\vac_{-1}\wedge u_2-q\vac_{-2} \wedge u_0 \wedge u_1
\]
is bar-invariant. Thus
\begin{eqnarray*}
G^{up}((1^2))=\langle (1^2)|, \quad G^{up}((2))=\langle(2)|-q\langle(1^2)|,
\end{eqnarray*}
and
\[
\langle(2)|=G^{up}((2))+qG^{up}((1^2)).
\]
Therefore the matrix $D=(d_{\lambda\mu})$ coincides with Example \ref{exalow}(2).\\
(3) \ By using the bar-invariance of $F_0,F_1$, we can obtain the bar-invariance of 
\[
F_0(\vac_{-2} \wedge u_0 \wedge u_1), \ F_1(\vac_{-2} \wedge u_0 \wedge u_1), \ F_0(\vac_{-1}\wedge u_2-q\vac_{-2} \wedge u_0 \wedge u_1).
\]
It is easy to compute $G^{up}$ and see that the matrix $D$ is coincide the Example\ref{exalow}(3).
\end{exa}
\subsection{LLT-Ariki type theorem on the $v$-Schur algebras}
\subsubsection{LLT-Ariki type theorem on the Hecke algebras}
\ We consider the Hecke algebra $H_n$ of type $A_n$. If $\zeta$ is not a root of unity, the simple modules of $H_n$ are indexed by the partitions of $n$. Let $S^\lambda$ be the simple module corresponding to a partition $\lambda$ of $n$. If $\zeta$ is an $\ell$-th root of unity, then the simple modules of $H_n$ are indexed by the $\ell$-regular partitions. Let $D^\mu$ be the simple module of $H_n$ at $\zeta=\sqrt[\ell]{1}$ corresponding to an $\ell$-regular partitions $\mu$ of $n$. The composition multiplicities
\[
d_{\lambda\mu}:=[S^\lambda,D^{\mu}]
\]
are called the decomposition numbers.\\
\indent Lascoux-Leclerc-Thibon conjectured in \cite{LLT} that this decomposition numbers are described by the global and crystal basis of integrable highest weight module $L(\Lambda_0)$ of $U_q(\widehat{\mf{sl}_\ell})$. More precisely, consider the submodule $L(\lambda_0)$ of $\mathcal{F}$ defined by
\[
L(\Lambda_0)=U_q^{-}(\widehat{\mf{sl}_\ell}))|\phi\rangle.
\] 
Then 
\[
\{G^{low}(\mu)|\mu \in \mathcal{P} \ \text{and $\mu$ is $\ell$-regular}\}.
\]
is a basis of $L(\lambda_0)$ and they have the expansion
\[
G^{low}(\mu)=|\mu\rangle+\sum_{\mu \triangleleft \lambda \in \mathcal{P}_n}d_{\lambda\mu}(q)|\lambda\rangle.
\]
They conjectured that the decomposition numbers 
\[
d_{\lambda'\mu'}=[S^{\lambda'}:D^{\mu'}]=d_{\lambda\mu}(1).
\]
This conjecture was proved by S. Ariki in \cite{Ari1}, and he proved the similar results on the cyclotomic Hecke algebras. For more detail, see also \cite{Ari2}.\\
\subsubsection{LLT-Ariki type theorem on the $v$-Schur algebras}
\ We consider the $v$-Schur algebra $\bb{S}(n)$. If $v$ is not a root of unity, the simple modules of $\bb{S}(n)$ are indexed by the partitions of $n$ and called the Weyl modules. Let $W^\lambda$ be the Weyl module corresponding to $\lambda \in \mathcal{P}_n$. If $v$ is an $\ell$-th root of unity, the Weyl modules $W^\lambda$ are generally reducible, but the simple modules of $\bb{S}(n)$ are indexed by the partitions of $n$. Let $L^\mu$ be these simple modules. The composition multiplicities 
\[
d_{\lambda\mu}=[W^\lambda:L^\mu]
\] 
are called the crystallized decomposition numbers.\\
\indent Varagnolo-Vasserot proved in \cite{VV1} that the crystallized decomposition numbers coincide $d_{\lambda\mu}(1)$. They extended the LLT-Ariki type theorem of the Hecke algebra to the $v$-Schur algebras.
\begin{thm}[Varagnolo-Vasserot\cite{VV1}]
Consider the Fock space $\mathcal{F}$ of $U_q(\widehat{\mf{sl}_\ell})$ and its lower global basis $\{G^{low}(\mu)|\mu \in \mathcal{P}\}$ and crystal basis $\{|\lambda\rangle|\lambda \in \mathcal{P}\}$. Let us consider the coefficients of 
\[
G^{low}(\mu)=|\mu\rangle+\sum_{\mu \triangleleft \lambda \in \mathcal{P}_n}d_{\lambda\mu}(q)|\lambda\rangle.
\]
Then we have
\[
d_{\lambda'\mu'}=[W^{\lambda'}:L^{\mu'}]=d_{\lambda\mu}(1).
\]
\end{thm}

\newpage
\section{Main Thorem}
\subsection{Main Theorem}

We will state the Main theorem on the coefficient $d_{(n),\mu}(q)$.\\
First, let us define the following partitions of $n$.
\begin{dfn}\label{defp}
Put $N=\left[\frac{n}{\ell}\right]$. The partitions $\mu_{i}^{(n)}$ is defined by the following,
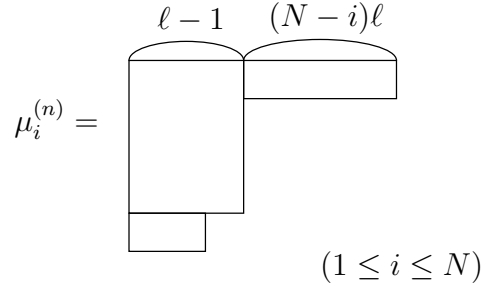
\[
\unitlength 0.1in
\begin{picture}(20.00,13.20)(14.00,-22.00)
%
\special{pn 8}%
\special{pa 2000 1200}%
\special{pa 2600 1200}%
\special{pa 2600 2000}%
\special{pa 2000 2000}%
\special{pa 2000 1200}%
\special{fp}%
%
\special{pn 8}%
\special{pa 2000 2000}%
\special{pa 2400 2000}%
\special{pa 2400 2200}%
\special{pa 2000 2200}%
\special{pa 2000 2000}%
\special{fp}%
%
\special{pn 8}%
\special{pa 2600 1200}%
\special{pa 3400 1200}%
\special{pa 3400 1400}%
\special{pa 2600 1400}%
\special{pa 2600 1200}%
\special{fp}%
%
\special{pn 8}%
\special{ar 2300 1200 300 100  3.1415927 6.2831853}%
%
\special{pn 8}%
\special{ar 3000 1200 400 110  3.1415927 6.2831853}%
\put(27.3000,-10.6000){\makebox(0,0)[lb]{$(N-i)\ell$}}%
\put(25.0000,-10.5000){\makebox(0,0)[rb]{$\ell-1$}}%
\put(14.0000,-14.0000){\makebox(0,0)[lt]{$\mu_i^{(n)}=$}}%
\put(30.0000,-22.0000){\makebox(0,0)[lt]{$(1 \le i \le N)$}}%
\end{picture}%
\]
and $\mu_0^{(n)}=(n)$.
\end{dfn}
\begin{rem}
These partitions of $n$ are constructed by removing a rim $\ell$-hook from the bottom to the first row. For example, let us consider the case $n=10, \ell=3$. Then $N=3$. In this case, $\mu_i^{(10)} \ (i=3,2,1,0)$ are the following,
\[
\mu_3^{(10)}=(2^5), \ 
\mu_2^{(10)}=(5,2^2,1), \ 
\mu_1^{(10)}=(8,2), \ 
\mu_0^{(10)}=(10).
\]
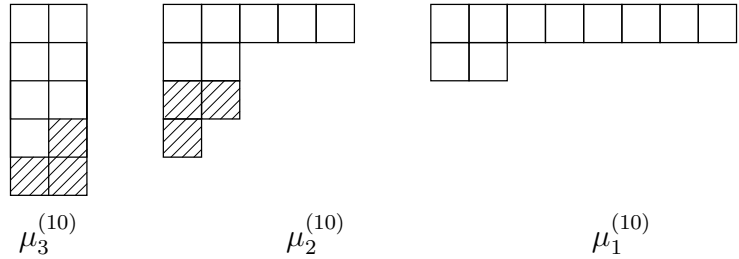
\[
\unitlength 0.1in
\begin{picture}(42.30,11.15)(3.70,-15.15)
%
\special{pn 8}%
\special{pa 800 400}%
\special{pa 1200 400}%
\special{pa 1200 1400}%
\special{pa 800 1400}%
\special{pa 800 400}%
\special{fp}%
%
\special{pn 8}%
\special{pa 1000 400}%
\special{pa 1000 1400}%
\special{fp}%
%
\special{pn 8}%
\special{pa 800 600}%
\special{pa 1200 600}%
\special{pa 1200 1200}%
\special{pa 800 1200}%
\special{pa 800 600}%
\special{fp}%
%
\special{pn 8}%
\special{pa 1200 1000}%
\special{pa 800 1000}%
\special{pa 800 800}%
\special{pa 1200 800}%
\special{pa 1200 1000}%
\special{fp}%
%
\special{pn 8}%
\special{pa 1600 400}%
\special{pa 1800 400}%
\special{pa 1800 600}%
\special{pa 1600 600}%
\special{pa 1600 400}%
\special{fp}%
%
\special{pn 8}%
\special{pa 1600 600}%
\special{pa 1800 600}%
\special{pa 1800 800}%
\special{pa 1600 800}%
\special{pa 1600 600}%
\special{fp}%
%
\special{pn 8}%
\special{pa 1600 800}%
\special{pa 1800 800}%
\special{pa 1800 1000}%
\special{pa 1600 1000}%
\special{pa 1600 800}%
\special{fp}%
%
\special{pn 8}%
\special{pa 1600 1000}%
\special{pa 1800 1000}%
\special{pa 1800 1200}%
\special{pa 1600 1200}%
\special{pa 1600 1000}%
\special{fp}%
%
\special{pn 8}%
\special{pa 1800 400}%
\special{pa 2000 400}%
\special{pa 2000 600}%
\special{pa 1800 600}%
\special{pa 1800 400}%
\special{fp}%
%
\special{pn 8}%
\special{pa 1800 600}%
\special{pa 2000 600}%
\special{pa 2000 800}%
\special{pa 1800 800}%
\special{pa 1800 600}%
\special{fp}%
%
\special{pn 8}%
\special{pa 1800 800}%
\special{pa 2000 800}%
\special{pa 2000 1000}%
\special{pa 1800 1000}%
\special{pa 1800 800}%
\special{fp}%
%
\special{pn 8}%
\special{pa 2000 400}%
\special{pa 2200 400}%
\special{pa 2200 600}%
\special{pa 2000 600}%
\special{pa 2000 400}%
\special{fp}%
%
\special{pn 8}%
\special{pa 2200 400}%
\special{pa 2400 400}%
\special{pa 2400 600}%
\special{pa 2200 600}%
\special{pa 2200 400}%
\special{fp}%
%
\special{pn 8}%
\special{pa 2400 400}%
\special{pa 2600 400}%
\special{pa 2600 600}%
\special{pa 2400 600}%
\special{pa 2400 400}%
\special{fp}%
%
\special{pn 8}%
\special{pa 3000 400}%
\special{pa 3200 400}%
\special{pa 3200 600}%
\special{pa 3000 600}%
\special{pa 3000 400}%
\special{fp}%
%
\special{pn 8}%
\special{pa 3000 600}%
\special{pa 3200 600}%
\special{pa 3200 800}%
\special{pa 3000 800}%
\special{pa 3000 600}%
\special{fp}%
%
\special{pn 8}%
\special{pa 4000 400}%
\special{pa 4200 400}%
\special{pa 4200 600}%
\special{pa 4000 600}%
\special{pa 4000 400}%
\special{fp}%
%
\special{pn 8}%
\special{pa 3200 400}%
\special{pa 3400 400}%
\special{pa 3400 600}%
\special{pa 3200 600}%
\special{pa 3200 400}%
\special{fp}%
%
\special{pn 8}%
\special{pa 3200 600}%
\special{pa 3400 600}%
\special{pa 3400 800}%
\special{pa 3200 800}%
\special{pa 3200 600}%
\special{fp}%
%
\special{pn 8}%
\special{pa 4200 400}%
\special{pa 4400 400}%
\special{pa 4400 600}%
\special{pa 4200 600}%
\special{pa 4200 400}%
\special{fp}%
%
\special{pn 8}%
\special{pa 3400 400}%
\special{pa 3600 400}%
\special{pa 3600 600}%
\special{pa 3400 600}%
\special{pa 3400 400}%
\special{fp}%
%
\special{pn 8}%
\special{pa 3600 400}%
\special{pa 3800 400}%
\special{pa 3800 600}%
\special{pa 3600 600}%
\special{pa 3600 400}%
\special{fp}%
%
\special{pn 8}%
\special{pa 3800 400}%
\special{pa 4000 400}%
\special{pa 4000 600}%
\special{pa 3800 600}%
\special{pa 3800 400}%
\special{fp}%
%
\special{pn 8}%
\special{pa 4400 400}%
\special{pa 4600 400}%
\special{pa 4600 600}%
\special{pa 4400 600}%
\special{pa 4400 400}%
\special{fp}%
\put(10.0000,-16.0000){\makebox(0,0){$\mu_3^{(10)}$}}%
\put(24.0000,-16.0000){\makebox(0,0){$\mu_2^{(10)}$}}%
\put(40.0000,-16.0000){\makebox(0,0){$\mu_1^{(10)}$}}%
%
\special{pn 4}%
\special{pa 980 1200}%
\special{pa 800 1380}%
\special{fp}%
\special{pa 1000 1240}%
\special{pa 840 1400}%
\special{fp}%
\special{pa 1000 1300}%
\special{pa 900 1400}%
\special{fp}%
\special{pa 1000 1360}%
\special{pa 960 1400}%
\special{fp}%
\special{pa 920 1200}%
\special{pa 800 1320}%
\special{fp}%
\special{pa 860 1200}%
\special{pa 800 1260}%
\special{fp}%
%
\special{pn 4}%
\special{pa 1200 1220}%
\special{pa 1020 1400}%
\special{fp}%
\special{pa 1160 1200}%
\special{pa 1000 1360}%
\special{fp}%
\special{pa 1100 1200}%
\special{pa 1000 1300}%
\special{fp}%
\special{pa 1040 1200}%
\special{pa 1000 1240}%
\special{fp}%
\special{pa 1200 1280}%
\special{pa 1080 1400}%
\special{fp}%
\special{pa 1200 1340}%
\special{pa 1140 1400}%
\special{fp}%
%
\special{pn 4}%
\special{pa 1200 1040}%
\special{pa 1040 1200}%
\special{fp}%
\special{pa 1180 1000}%
\special{pa 1000 1180}%
\special{fp}%
\special{pa 1120 1000}%
\special{pa 1000 1120}%
\special{fp}%
\special{pa 1060 1000}%
\special{pa 1000 1060}%
\special{fp}%
\special{pa 1200 1100}%
\special{pa 1100 1200}%
\special{fp}%
\special{pa 1200 1160}%
\special{pa 1160 1200}%
\special{fp}%
%
\special{pn 4}%
\special{pa 1740 800}%
\special{pa 1600 940}%
\special{fp}%
\special{pa 1790 810}%
\special{pa 1610 990}%
\special{fp}%
\special{pa 1800 860}%
\special{pa 1660 1000}%
\special{fp}%
\special{pa 1800 920}%
\special{pa 1720 1000}%
\special{fp}%
\special{pa 1680 800}%
\special{pa 1600 880}%
\special{fp}%
%
\special{pn 4}%
\special{pa 2000 840}%
\special{pa 1840 1000}%
\special{fp}%
\special{pa 1980 800}%
\special{pa 1800 980}%
\special{fp}%
\special{pa 1920 800}%
\special{pa 1800 920}%
\special{fp}%
\special{pa 1860 800}%
\special{pa 1800 860}%
\special{fp}%
\special{pa 2000 900}%
\special{pa 1900 1000}%
\special{fp}%
\special{pa 2000 960}%
\special{pa 1960 1000}%
\special{fp}%
%
\special{pn 4}%
\special{pa 1780 1000}%
\special{pa 1600 1180}%
\special{fp}%
\special{pa 1800 1040}%
\special{pa 1640 1200}%
\special{fp}%
\special{pa 1800 1100}%
\special{pa 1700 1200}%
\special{fp}%
\special{pa 1800 1160}%
\special{pa 1760 1200}%
\special{fp}%
\special{pa 1720 1000}%
\special{pa 1600 1120}%
\special{fp}%
\special{pa 1660 1000}%
\special{pa 1600 1060}%
\special{fp}%
\end{picture}%
\]
\end{rem}
\quad {} \\
\begin{thm}
The crystal basis element $\langle(n)|$ is expanded by $G^{up}$ as the following;
\[
\langle(n)|=\sum_{i=0}^{N}q^iG^{up}(\mu_i^{(n)}),
\]
i.e. the $q$-decomposition numbers
\[
d_{(n),\mu}(q)=\left\{
\begin{array}{ll}
q^i & \text{if} \ \mu=\mu_{i}^{(n)} \\
0 & \text{otherwise}
\end{array}
\right..
\]
\begin{cor}
The decomposition numbers 
\[
d_{(n),\mu}=[W_{(1^n)}:L_\mu]=\left\{
\begin{array}{ll}
1 & \text{if} \ \mu=\mu_{i}^{(n)} \\
0 & \text{otherwise}
\end{array}
\right..
\]
\end{cor}
\end{thm}
\subsection{Proof of Main Theorem}
\subsubsection{Key Lemmas}
\indent First, we can describe the action of $E_i$ on the upper global base $\{G^{up}\}$. Let us set 
\[
\varepsilon_i(\mu)=\max\{k \ge 0|\widehat{e_i}^k\mu \neq 0\}.
\]
\begin{lem}[Kashiwara \cite{Kas2}]\label{l1} The element $E_i$ in $U_q(\widehat{sl}_\ell)$ acts on $G^{up}$ as the following,
\[
E_iG^{up}(\mu)=[\varepsilon_i(\mu)]G^{up}(\widetilde{e_i}\mu)+\sum_{\varepsilon_i(\nu)<\varepsilon_i(\mu)-1}b_{\mu\nu}^iG^{up}(\nu).
\]
Especially, if $\varepsilon_i(\mu)=1$, then 
\[
E_iG^{up}(\mu)=G^{up}(\widetilde{e_i}\mu).
\]
\end{lem}
\indent Secondly, we have the following lemma on a property of $\displaystyle \bigcap_j\Ker(E_j)$.
\begin{lem}\label{l2}
For any $\displaystyle x \in \bigcap_{j}\Ker(E_j) \subset \mathcal{F}^\vee$, we have the following expansion;
\[
x=\sum_{\lambda \in \mathcal{P}}b_{x,\lambda}G^{up}(\ell\lambda).
\]
\end{lem}
\indent Thirdly, we consider the coefficients of $G^{up}(\ell\lambda)$ in the expansion of $\langle (n)|$.
\begin{lem}\label{l3} \quad {} \\
(1) (Kashiwara \cite{Kas}) \ The following expansions hold,
\begin{eqnarray*}
&{}&G^{low}(\vac_{-m-1} \wedge u_{\ell\lambda_m-m} \wedge \cdots \wedge u_{\ell\lambda_1-1} \wedge u_{\ell\lambda_0}) \\
&=&\sum{a_{j_m,j_{m-1}, \cdots ,j_0}(q)}\vac_{-m-1} \wedge u_{j_m+\ell\lambda_m-{\ell}m} \wedge u_{j_{m-1}+\ell\lambda_{m-1}-\ell(m-1)} \wedge \cdots \wedge u_{j_0+\ell\lambda_0}
\end{eqnarray*}
where the sum runs on the index $(j_m,j_{m-1}, \cdots ,{j_0})$ satisfied
\[
(0,\ell-1,2(\ell-1), \cdots ,m(\ell-1)) \le (j_m,j_{m-1}, \cdots ,j_0) \le (m(\ell-1),(m-1)(\ell-1), \cdots ,0)
\]
and
\[
(-m,-m+1, \cdots ,0) \le (j_m+\ell\lambda_m-{\ell}m,j_{m-1}+\ell\lambda_{m-1}-\ell(m-1), \cdots ,j_0+\ell\lambda_0).
\]
(2) \ The element $|(n)\rangle$ does not appear in the expansion of $G^{low}(\ell\lambda) \ (\lambda \in \mathcal{P})$ with respect to the crystal base. \\
(3) \ The elements $G^{up}(\ell\lambda) \ (\lambda \in \mathcal{P})$ do not appear in the expansion of $\langle(n)|$ with respect to the upper global base.
\end{lem}

\indent Forthly, we can describe the action of the Kashiwara operators $\widetilde{e_j}$ on $\mu_i^{(n)}$ by using the Misra-Miwa's Theorem \ref{MM}.
\begin{lem}\label{l4}
The action of $\widetilde{e_j}$ on $\mu_i^{(n)}$ is obtained by the following; \\
(1) \ If $n \not\equiv 0 \ (\mod{\ell})$, then
\[
\widetilde{e_j}(\mu_i^{(n)})=\left\{
\begin{array}{ll}
0 & \text{if} \ j \not\equiv n-1 \\
\mu_i^{(n-1)} & \text{if} \ j \equiv n-1
\end{array}
\right.
\]
for $0 \le i \le N$. \\
(2) \ If $n \equiv 0 \ (\mod{\ell})$, then
\[
\widetilde{e_j}(\mu_i^{(n)})=\left\{
\begin{array}{ll}
0 & \text{if} \ j \not\equiv n-1 \\
\mu_{i-1}^{(n-1)} & \text{if} \ j \equiv n-1
\end{array}
\right.
\]
for $1 \le i \le N$. And $\widetilde{e_j}\mu_0^{(n)}=0$ for any $0 \le j \le \ell-1$. \\
(3) \ Especially, 
\begin{eqnarray*}
\varepsilon_j(\mu_i^{(n)})=\left\{
\begin{array}{ll}
1 & \text{if} \ j \equiv n-1 \\
0 & \text{otherwise}  
\end{array}
\right.
\end{eqnarray*}
for $1 \le i \le N$. And 
\[
\varepsilon_j(\mu_0^{(n)})=\left\{
\begin{array}{ll}
1 & \text{if} \ j \equiv n-1 \ \text{and} \ n \not\equiv 0 \\
0 & \text{otherwise}  
\end{array}
\right..
\]
\end{lem}
\begin{exa}
We consider the case $n=10, \ell=3$. Then
\begin{eqnarray*}
&{}& \widetilde{e_0}\mu_2^{(10)}=(5,2,2)=\mu_2^{(9)}, \\
&{}& \widetilde{e_0}\mu_2^{(9)}=(5,2,1)=\mu_1^{(8)}.
\end{eqnarray*}
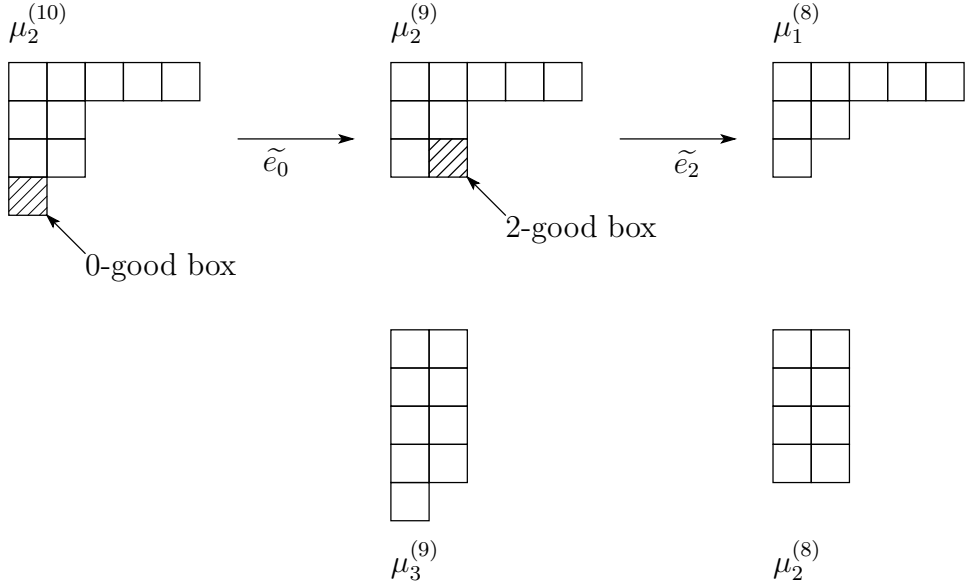
\[
\unitlength 0.1in
\begin{picture}(50.00,27.70)(8.00,-29.00)
%
\special{pn 8}%
\special{pa 800 400}%
\special{pa 1000 400}%
\special{pa 1000 600}%
\special{pa 800 600}%
\special{pa 800 400}%
\special{fp}%
%
\special{pn 8}%
\special{pa 800 600}%
\special{pa 1000 600}%
\special{pa 1000 800}%
\special{pa 800 800}%
\special{pa 800 600}%
\special{fp}%
%
\special{pn 8}%
\special{pa 800 800}%
\special{pa 1000 800}%
\special{pa 1000 1000}%
\special{pa 800 1000}%
\special{pa 800 800}%
\special{fp}%
%
\special{pn 8}%
\special{pa 800 1000}%
\special{pa 1000 1000}%
\special{pa 1000 1200}%
\special{pa 800 1200}%
\special{pa 800 1000}%
\special{fp}%
%
\special{pn 8}%
\special{pa 1000 400}%
\special{pa 1200 400}%
\special{pa 1200 600}%
\special{pa 1000 600}%
\special{pa 1000 400}%
\special{fp}%
%
\special{pn 8}%
\special{pa 1000 600}%
\special{pa 1200 600}%
\special{pa 1200 800}%
\special{pa 1000 800}%
\special{pa 1000 600}%
\special{fp}%
%
\special{pn 8}%
\special{pa 1000 800}%
\special{pa 1200 800}%
\special{pa 1200 1000}%
\special{pa 1000 1000}%
\special{pa 1000 800}%
\special{fp}%
%
\special{pn 8}%
\special{pa 1200 400}%
\special{pa 1400 400}%
\special{pa 1400 600}%
\special{pa 1200 600}%
\special{pa 1200 400}%
\special{fp}%
%
\special{pn 8}%
\special{pa 1400 400}%
\special{pa 1600 400}%
\special{pa 1600 600}%
\special{pa 1400 600}%
\special{pa 1400 400}%
\special{fp}%
%
\special{pn 8}%
\special{pa 1600 400}%
\special{pa 1800 400}%
\special{pa 1800 600}%
\special{pa 1600 600}%
\special{pa 1600 400}%
\special{fp}%
\put(8.0000,-3.0000){\makebox(0,0)[lb]{$\mu_2^{(10)}$}}%
%
\special{pn 4}%
\special{pa 1000 1040}%
\special{pa 840 1200}%
\special{fp}%
\special{pa 980 1000}%
\special{pa 800 1180}%
\special{fp}%
\special{pa 920 1000}%
\special{pa 800 1120}%
\special{fp}%
\special{pa 860 1000}%
\special{pa 800 1060}%
\special{fp}%
\special{pa 1000 1100}%
\special{pa 900 1200}%
\special{fp}%
\special{pa 1000 1160}%
\special{pa 960 1200}%
\special{fp}%
%
\special{pn 8}%
\special{pa 1200 1400}%
\special{pa 1000 1200}%
\special{fp}%
\special{sh 1}%
\special{pa 1000 1200}%
\special{pa 1033 1261}%
\special{pa 1038 1238}%
\special{pa 1061 1233}%
\special{pa 1000 1200}%
\special{fp}%
\put(12.0000,-14.0000){\makebox(0,0)[lt]{$0$-good box}}%
%
\special{pn 8}%
\special{pa 2800 400}%
\special{pa 3000 400}%
\special{pa 3000 600}%
\special{pa 2800 600}%
\special{pa 2800 400}%
\special{fp}%
%
\special{pn 8}%
\special{pa 2800 600}%
\special{pa 3000 600}%
\special{pa 3000 800}%
\special{pa 2800 800}%
\special{pa 2800 600}%
\special{fp}%
%
\special{pn 8}%
\special{pa 2800 800}%
\special{pa 3000 800}%
\special{pa 3000 1000}%
\special{pa 2800 1000}%
\special{pa 2800 800}%
\special{fp}%
%
\special{pn 8}%
\special{pa 3000 400}%
\special{pa 3200 400}%
\special{pa 3200 600}%
\special{pa 3000 600}%
\special{pa 3000 400}%
\special{fp}%
%
\special{pn 8}%
\special{pa 3000 600}%
\special{pa 3200 600}%
\special{pa 3200 800}%
\special{pa 3000 800}%
\special{pa 3000 600}%
\special{fp}%
%
\special{pn 8}%
\special{pa 3000 800}%
\special{pa 3200 800}%
\special{pa 3200 1000}%
\special{pa 3000 1000}%
\special{pa 3000 800}%
\special{fp}%
%
\special{pn 8}%
\special{pa 3200 400}%
\special{pa 3400 400}%
\special{pa 3400 600}%
\special{pa 3200 600}%
\special{pa 3200 400}%
\special{fp}%
%
\special{pn 8}%
\special{pa 3400 400}%
\special{pa 3600 400}%
\special{pa 3600 600}%
\special{pa 3400 600}%
\special{pa 3400 400}%
\special{fp}%
%
\special{pn 8}%
\special{pa 3600 400}%
\special{pa 3800 400}%
\special{pa 3800 600}%
\special{pa 3600 600}%
\special{pa 3600 400}%
\special{fp}%
%
\special{pn 8}%
\special{pa 3180 800}%
\special{pa 3000 980}%
\special{fp}%
\special{pa 3200 840}%
\special{pa 3040 1000}%
\special{fp}%
\special{pa 3200 900}%
\special{pa 3100 1000}%
\special{fp}%
\special{pa 3200 960}%
\special{pa 3160 1000}%
\special{fp}%
\special{pa 3120 800}%
\special{pa 3000 920}%
\special{fp}%
\special{pa 3060 800}%
\special{pa 3000 860}%
\special{fp}%
%
\special{pn 8}%
\special{pa 3400 1200}%
\special{pa 3200 1000}%
\special{fp}%
\special{sh 1}%
\special{pa 3200 1000}%
\special{pa 3233 1061}%
\special{pa 3238 1038}%
\special{pa 3261 1033}%
\special{pa 3200 1000}%
\special{fp}%
\put(34.0000,-12.0000){\makebox(0,0)[lt]{$2$-good box}}%
\put(28.0000,-3.0000){\makebox(0,0)[lb]{$\mu_2^{(9)}$}}%
\put(48.0000,-3.0000){\makebox(0,0)[lb]{$\mu_1^{(8)}$}}%
%
\special{pn 8}%
\special{pa 4800 400}%
\special{pa 5000 400}%
\special{pa 5000 600}%
\special{pa 4800 600}%
\special{pa 4800 400}%
\special{fp}%
%
\special{pn 8}%
\special{pa 4800 600}%
\special{pa 5000 600}%
\special{pa 5000 800}%
\special{pa 4800 800}%
\special{pa 4800 600}%
\special{fp}%
%
\special{pn 8}%
\special{pa 5000 400}%
\special{pa 5200 400}%
\special{pa 5200 600}%
\special{pa 5000 600}%
\special{pa 5000 400}%
\special{fp}%
%
\special{pn 8}%
\special{pa 5000 600}%
\special{pa 5200 600}%
\special{pa 5200 800}%
\special{pa 5000 800}%
\special{pa 5000 600}%
\special{fp}%
%
\special{pn 8}%
\special{pa 5200 400}%
\special{pa 5400 400}%
\special{pa 5400 600}%
\special{pa 5200 600}%
\special{pa 5200 400}%
\special{fp}%
%
\special{pn 8}%
\special{pa 5400 400}%
\special{pa 5600 400}%
\special{pa 5600 600}%
\special{pa 5400 600}%
\special{pa 5400 400}%
\special{fp}%
%
\special{pn 8}%
\special{pa 5600 400}%
\special{pa 5800 400}%
\special{pa 5800 600}%
\special{pa 5600 600}%
\special{pa 5600 400}%
\special{fp}%
%
\special{pn 8}%
\special{pa 4800 800}%
\special{pa 5000 800}%
\special{pa 5000 1000}%
\special{pa 4800 1000}%
\special{pa 4800 800}%
\special{fp}%
%
\special{pn 8}%
\special{pa 2000 800}%
\special{pa 2600 800}%
\special{fp}%
\special{sh 1}%
\special{pa 2600 800}%
\special{pa 2533 780}%
\special{pa 2547 800}%
\special{pa 2533 820}%
\special{pa 2600 800}%
\special{fp}%
%
\special{pn 8}%
\special{pa 4000 800}%
\special{pa 4600 800}%
\special{fp}%
\special{sh 1}%
\special{pa 4600 800}%
\special{pa 4533 780}%
\special{pa 4547 800}%
\special{pa 4533 820}%
\special{pa 4600 800}%
\special{fp}%
\put(22.0000,-9.2000){\makebox(0,0){$\widetilde{e_0}$}}%
\put(43.5000,-9.3000){\makebox(0,0){$\widetilde{e_2}$}}%
%
\special{pn 8}%
\special{pa 2800 1800}%
\special{pa 3000 1800}%
\special{pa 3000 2000}%
\special{pa 2800 2000}%
\special{pa 2800 1800}%
\special{fp}%
%
\special{pn 8}%
\special{pa 2800 2000}%
\special{pa 3000 2000}%
\special{pa 3000 2200}%
\special{pa 2800 2200}%
\special{pa 2800 2000}%
\special{fp}%
%
\special{pn 8}%
\special{pa 3000 1800}%
\special{pa 3200 1800}%
\special{pa 3200 2000}%
\special{pa 3000 2000}%
\special{pa 3000 1800}%
\special{fp}%
%
\special{pn 8}%
\special{pa 3000 2000}%
\special{pa 3200 2000}%
\special{pa 3200 2200}%
\special{pa 3000 2200}%
\special{pa 3000 2000}%
\special{fp}%
%
\special{pn 8}%
\special{pa 2800 2400}%
\special{pa 3000 2400}%
\special{pa 3000 2600}%
\special{pa 2800 2600}%
\special{pa 2800 2400}%
\special{fp}%
%
\special{pn 8}%
\special{pa 2800 2200}%
\special{pa 3000 2200}%
\special{pa 3000 2400}%
\special{pa 2800 2400}%
\special{pa 2800 2200}%
\special{fp}%
%
\special{pn 8}%
\special{pa 2800 2600}%
\special{pa 3000 2600}%
\special{pa 3000 2800}%
\special{pa 2800 2800}%
\special{pa 2800 2600}%
\special{fp}%
%
\special{pn 8}%
\special{pa 3000 2200}%
\special{pa 3200 2200}%
\special{pa 3200 2400}%
\special{pa 3000 2400}%
\special{pa 3000 2200}%
\special{fp}%
%
\special{pn 8}%
\special{pa 3000 2400}%
\special{pa 3200 2400}%
\special{pa 3200 2600}%
\special{pa 3000 2600}%
\special{pa 3000 2400}%
\special{fp}%
%
\special{pn 8}%
\special{pa 4800 1800}%
\special{pa 5000 1800}%
\special{pa 5000 2000}%
\special{pa 4800 2000}%
\special{pa 4800 1800}%
\special{fp}%
%
\special{pn 8}%
\special{pa 4800 2000}%
\special{pa 5000 2000}%
\special{pa 5000 2200}%
\special{pa 4800 2200}%
\special{pa 4800 2000}%
\special{fp}%
%
\special{pn 8}%
\special{pa 5000 1800}%
\special{pa 5200 1800}%
\special{pa 5200 2000}%
\special{pa 5000 2000}%
\special{pa 5000 1800}%
\special{fp}%
%
\special{pn 8}%
\special{pa 5000 2000}%
\special{pa 5200 2000}%
\special{pa 5200 2200}%
\special{pa 5000 2200}%
\special{pa 5000 2000}%
\special{fp}%
%
\special{pn 8}%
\special{pa 4800 2200}%
\special{pa 5000 2200}%
\special{pa 5000 2400}%
\special{pa 4800 2400}%
\special{pa 4800 2200}%
\special{fp}%
%
\special{pn 8}%
\special{pa 5000 2200}%
\special{pa 5200 2200}%
\special{pa 5200 2400}%
\special{pa 5000 2400}%
\special{pa 5000 2200}%
\special{fp}%
%
\special{pn 8}%
\special{pa 4800 2400}%
\special{pa 5000 2400}%
\special{pa 5000 2600}%
\special{pa 4800 2600}%
\special{pa 4800 2400}%
\special{fp}%
%
\special{pn 8}%
\special{pa 5000 2400}%
\special{pa 5200 2400}%
\special{pa 5200 2600}%
\special{pa 5000 2600}%
\special{pa 5000 2400}%
\special{fp}%
\put(28.0000,-29.0000){\makebox(0,0)[lt]{$\mu_3^{(9)}$}}%
\put(48.0000,-29.0000){\makebox(0,0)[lt]{$\mu_2^{(8)}$}}%
\end{picture}%
\]
\end{exa}
\subsubsection{Proof of Main Theorem}
\indent First, we show that 
\[
D_n:=\vac_{-1} \wedge u_n-\sum_{i=0}^{N}q^iG^{up}(\mu_i^{(n)}) \in \bigcap_{j=0}^{\ell-1}\Ker(E_j)
\]
by the induction with respect to $n$.\\
\indent Let us assume that $D_{n-1} \in \cap\Ker(E_j)$. \\
\indent If $n \equiv 0 \ (\mod{\ell})$, then
\[
E_j\vac_{-1} \wedge u_n=0 \ (0 \le j \le \ell-2)
\]
and 
\[
E_jG^{up}(\mu_i^{(n)})=0 \ (0 \le j \le \ell-2)
\]
by Lemma \ref{l1} and Lemma \ref{l4}. Therefore, $E_j(D_n)=0$ for $0 \le j \le \ell-1$. On the other hand, 
\[
E_{\ell-1}(D_n)=q\vac_{-1} \wedge u_{n-1}-\sum_{i=1}^{N}q^iG^{up}(u_{i-1}^{(n-1)})=qD_{n-1} \in \bigcap_{j=1}^{\ell-1}\Ker(E_j) 
\]
by the induction hypothesis. Thus $\displaystyle D_n \in \bigcap_{j=1}^{\ell-1}\Ker(E_j)$. \\
\indent If $n \not\equiv 0 \ (\mod{\ell})$, then similarly as above
\begin{eqnarray*}
E_j(D_n)=0
\end{eqnarray*}
for $j \not\equiv n-1$, and
\begin{eqnarray*}
E_j(D_n)&=&\vac_{-1} \wedge u_{n-1}-\sum_{i=0}^{N}q^iG^{up}(\mu_i^{(n-1)})=D_{n-1} \in \bigcap_{j=1}^{\ell-1}\Ker(E_j)
\end{eqnarray*}
for $j \equiv n-1$. Thus $\displaystyle D_n \in \bigcap_{j=1}^{\ell-1}\Ker(E_j)$. \\
\indent By Lemma \ref{l2}, $D_n$ is expanded by $\{G^{up}(\ell\lambda)|\lambda \in \mathcal{P}\}$. But this is contradiction except for $D_n=0$ by Lemma \ref{l3}. Thus the proof of theorem is complete.

\subsubsection{Two Remarks}
\begin{rem}[Prior Results and Conjectures by H. Miyachi.]
The results
\[
d_{(n),\mu}(1)=[W^{(1^n)}:L^\mu]=\left\{
\begin{array}{ll}
1 & \text{if} \ \mu=\mu_{i}^{(n)} \\
0 & \text{otherwise}
\end{array}
\right.
\]
were proved in \cite[Lemma 12.2.4 and Corollary 12.2.6]{Mi}. The $q$-decomposition numbers $d_{(n)\mu}(q)$ were also conjectured in \cite[Conjecture 12.2.19]{Mi}, and are calculated in this paper.
\end{rem}
\begin{rem}[Combinatorics of decomposition numbers]
We can prove only
\[
d_{(n),\mu}(1)=1 \ \text{if $\mu=\mu_i^{(n)}$}
\]
by using the combinatorial techniques of the decomposition numbers. But I cannot prove $d_{(n),\mu}=0$ unless $\mu=\mu_i^{(n)}$ by using only combinatorics. \\
\indent In this remark, we will show the sketch of the above combinatorial proof. We use the three combinatorial Lemmas about decomposition numbers. \\
\indent First, we can reduce the decomposition numbers $d_{\lambda\mu}$ of the $v$-Schur algebras to the decomposition numbers $d_{\tilde{\lambda}\tilde{\mu}}$ of the Hecke algebras. But the size of partitions are quite large. The following Lemma is the special case $\lambda=(n)$.
\begin{lem}[Leclerc \cite{L}]\label{LL1}
Consider the $v$-Schur algebra at $v \in \sqrt[\ell]{1}$ and the decomposition number $d_{(n),\mu}$ such that $\mu$ is not $\ell$-regular and $\mu$ has $m$ rows. Then we have a unique decomposition $\mu=\mu^{(1)}+\ell\mu^{(0)}$ such that $\mu^{(1)}$ is an $\ell$-restricted partition. Let $\tilde{\mu}$ be the partition 
\[
\widetilde{\mu}=((\ell-1)(m-1),(\ell-1)(m-2), \cdots ,\ell-1,0)+w_0(\mu^{(1)})+\ell\mu^{(0)},
\]
where $w_0(\mu^{(1)})=(\mu_m,\mu_{m-1}, \cdots ,\mu_1)$ for $\mu^{(1)}=(\mu_1, \cdots ,\mu_{m-1},\mu_m)$.
\indent Let $\widetilde{(n)}$ be the partition
\[
(n+(\ell-1)(m-1),(\ell-1)(m-1), \cdots ,(\ell-1)(m-1)).
\]
Then
\[
d_{(n),\mu}=d_{\widetilde{(n)},\widetilde{\mu}}.
\]
\end{lem}
\indent Secondly, the following lemma is a special case of the row and column removal formula.  See \cite[6.4 Rule 8]{Mat}.
\begin{lem}[Row and Column removal formula]\label{LL2}
Let $\lambda,\mu$ be two partitions and $\lambda',\mu'$ be their conjugate partitions.\\
(1) \ If $\lambda_i=\mu_i \ (1 \le i \le r)$, then $d_{\lambda\mu}=d_{(\lambda_{r+1},\lambda{r+2}, \cdots),(\mu_{r+1},\mu_{r+2}, \cdots)}$. \\
(2) \ If ${\lambda'}_i={\mu'}_i \ (1 \le i \le r)$, then $d_{\lambda'\mu'}=d_{({\lambda'}_{r+1},{\lambda'}_{r+2}, \cdots),({\mu'}_{r+1},{\mu'}_{r+2}, \cdots)}$.
\end{lem}
Thirdly, the following lemma describe the relationship between the Kleshchev-Mullineux involution and the decomposition numbers of the Hecke algebras. See \cite[6.4 Rule 11]{Mat}.
\begin{lem}[Kleshchev-Mullineux involution]\label{LL3}
Suppose the partition $\nu=\widetilde{f_{i_s}} \cdots \widetilde{f_{i_2}}\widetilde{f_{i_1}} \phi$ as the Misra-Miwa's theorem \ref{MM}. Let us define the partition
\[
\mathbf{m}(\nu)=\widetilde{f_{-i_s}} \cdots \widetilde{f_{-i_2}}\widetilde{f_{-i_1}} \phi
\]
where the indeices are to be read modulo $\ell$. Then
\[
d_{\lambda\nu}=d_{\lambda'\mathbf{m}(\nu)}.
\]
\end{lem}
\indent We sketch the proof and give one example. \\
\indent First, by using the Lemma \ref{LL1}, we can obtain $d_{(n),\mu_i^{(n)}}=d_{\widetilde{(n)},\widetilde{\mu_i^{(n)}}}$. Next we can cut off the several rows and columns from the two partitions by Lemma \ref{LL2}. And by Lemma \ref{LL3}, we can obtain the new partitions. In this case, it is easy to describe the image of Kleshchev-Mullineux involution. And we can cut off more several rows and columns. By repeating this step, two partitions can be coincide. Thus the decomposition numbers are equal to $1$.
\begin{exa}
We set $\ell=3, n=12$. We consider the decomposition number $d_{(12),\mu}$ where 
\[
\mu=(5,2,2,2,1).
\]
Then by the Lemma \ref{LL1}, we have
\[
d_{(12),(5,2,2,2,1)}=d_{(20,8,8,8,8),(20,14,10,6,2)}.
\]
Next, by using the Lemma \ref{LL2}, we have
\[
d_{(20,8,8,8,8),(20,14,10,6,2)}=d_{(6,6,6,6),(12,8,4)}.
\]
Note that $\mathbf{m}((12,8,4))=(6,6,4,4,2,2)$. Thus 
\[
d_{(6,6,6,6),(12,8,4)}=d_{(4,4,4,4,4,4),(6,6,4,4,2,2)}
\]
by using the Lemma \ref{LL3}. By using the Lemma \ref{LL2} and the Lemma \ref{LL3} again and again, we have
\[
d_{(4,4,4,4,4,4),(6,6,4,4,2,2)}=d_{(6,6),(8,4)}=d_{(2,2),(4)}=d_{(2,2),(2,2)}=1.
\]
\end{exa}

\end{rem}


\newpage
\begin{thebibliography}{GGOR}
\bibitem[Ari1]{Ari1}S. Ariki, \textit{On the decomposition numbers of the
Hecke algebra of $G(m,1,n)$},  J. Math. Kyoto Univ.  36  (1996),  no. 4,
789--808
\bibitem[Ari2]{Ari2}S. Ariki, \textit{Representations of quantum algebras
and combinatorics of Young tableaux}, University Lecture Series, 26.
American Mathematical Society, Providence, RI, 2002

\bibitem[AST]{AST}T. Arakawa, T. Suzuki, A. Tsuchiya, \textit{Degenerate
double affine Hecke algebra and conformal field theory}, Topological field
theory, primitive forms and related topics (Kyoto, 1996),  1--34, Progr.
Math., 160, Birkh\"{a}user Boston, Boston, MA, 1998, q-alg/9710031

\bibitem[BLM]{BLM}A. Beilinson, G. Lusztig, R. MacPherson, \textit{A geometric setting for the quantum deformation of $GL_n$}, Duke math. J., 61, 1990, 655-677

\bibitem[CG]{CG}N. Chriss, V. Ginzburg, \textit{Representation Theory and
Complex Geometry}, Birkh\"{a}user, 1997


\bibitem[Ch1]{Ch1}I. Cherednik, \textit{Double affine Hecke algebras,
Knizhnik-Zamolodchikov equations, and Macdonald's operators}, IMRS 1992,
no.9. 171-180
\bibitem[Ch2]{Ch2}I. Cherednik, \textit{Double affine Hecke algebras and
Macdonald conjectures}, Annals of math, 141, 1995
\bibitem[Ch3]{Ch3}I. Cherednik, \textit{Double Affine Hecke Algebras},
London Mathematical Society Lecture Note Series, 2005
\bibitem[Ch4]{Ch4}I. Cherednik, \textit{Elliptic quantum many-body problem
and double affine Knizhnik-Zamolodchikov equation},  Comm. Math. Phys.
169  (1995), no.2, 441--461

\bibitem[DJ]{DJ}R. Dipper, G. James, \textit{The $q$-Schur algebra}, Proc. AMS (3), 59, 1989, 23-50

\bibitem[EG]{EG}P. Etingof, V. Ginzburg, \textit{Symplectic reflection
algebras, Calogero-Moser space, and deformed Harish-Chandra homomorphism},
Invent. Math., 147, 2002, no.2, 243--348

\bibitem[G]{G}V. Ginzburg, \textit{Geometric methods in representation
theory of Hecke algebras and quantum groups}, in "Representation theories
and geometry" (Montreal, PQ, 1997), 127-183, Kluwer Acad. 1998
(math.AG/9802004)

\bibitem[GGOR]{GGOR} V. Ginzburg, N. Guay, E. Opdam, R. Rouquier,
\textit{On the category $O$ for rational Cherednik algebras},  Invent.
Math.  154  (2003),  no. 3, 617--651, math.RT/0212036

\bibitem[GL]{GL}J. Graham, G. Lehrer, \textit{Cellular algebras}, Invent. Math., 123, 1996, 1-34

\bibitem[GSI]{GS1}I. Gordon, J. T. Stafford, \textit{Rational Cherednik
algebras and Hilbert schemes I}, math.RT/0407516
\bibitem[GSII]{GS2}I. Gordon, J. T. Stafford, \textit{Rational Cherednik
algebras and Hilbert schemes II, Representations and Sheaves},
math.RA/0410293

\bibitem[Ha]{Ha}T. Hayashi, \textit{$q$-analogues of Clifford and Weyl algebras---spinor and oscillator representations of quantum enveloping algebras},  Comm. Math. Phys., 127, 1990, no.1, 129--144

\bibitem[IM]{IM}N. Iwahori, H. Matsumoto, \textit{On some Bruhat
decomposition and the structure of the Hecke rings of $p$-adic Chevalley
groups}, Publications Mathematiques de l'IHES, 25 (1965), p. 5-48

\bibitem[Kas1]{Kas} M. Kashiwara, \textit{On Level Zero Representations of
Quantized Affine Algebras},  Duke Math. J.  112  (2002),  no. 1, 117--175,
math.QA/0010293
\bibitem[Kas2]{Kas2} M. Kashiwara, \textit{Global crystal bases of quantum groups},  Duke Math. J.  69  (1993),  no. 2, 455--485

\bibitem[Kasa]{Kasa}M. Kasatani, \textit{Subrepresentations in the
Polynomial Representation of the Double Affine Hecke Algebra of type $GL_n$
at $t^{k+1}q^{r-1}=1$},  Int. Math. Res. Not.  2005,  no. 28, 1717--1742,
math.QA/0501272

\bibitem[KL]{KL}D. Kazhdan, G. Lusztig, \textit{Proof of Delingre-Langlands
Conjecture for Hecke Algebras}, Invent. math. \textbf{87}, 153-215, 1987

\bibitem[KMS]{KMS}M. Kashiwara, T. Miwa, E. Stern, \textit{Decomposition
of $q$-deformed Fock spaces},  Selecta Math. (N.S.)  1  (1995),  no. 4,
787--805, q-alg/9508006

\bibitem[L]{L}B. Leclerc, \textit{Decomposition numbers and canonical
bases}, math.QA/9902006

\bibitem[LLT]{LLT}A. Lascoux, B. Leclerc, J. Thibon, \textit{Hecke
algebras at roots of unity and crystal bases of quantum affine algebras},
Comm. Math. Phys.  181  (1996),  no. 1, 205--263

\bibitem[LT]{LT}B. Leclerc, J. Thibon, \textit{Canonical bases of
$q$-deformed Fock spaces},  Internat. Math. Res. Notices  1996,  no. 9,
447--456, q-alg/9602025

\bibitem[Lus]{Lus}G. Lusztig, \textit{Affine Hecke algebras and their
graded version},  J. Amer. Math. Soc.  2  (1989),  no. 3, 599--635

\bibitem[Mat]{Mat}A. Mathas, \textit{Iwahori-Hecke algebras and Schur
algebras of the symmetric groups}, University Lecture Series, 15. American
Mathematical Society, Providence, RI, 1999

\bibitem[Mi]{Mi}H. Miyachi, \textit{Unipotent Blocks of Finite General Linear Groups in Non-defining Characteristic}

\bibitem[MM]{MM}K. C. Misra, T. Miwa, \textit{Crystal base for the basic representation of $U_q(\widehat{sl}_n)$}, Comm. Math. Phys. 134, 79-88, 1990

\bibitem[Ra1]{R1}A. Ram, \textit{Affine Hecke Algebras and generalized
standard Young tableaux}, J. \ Algebra, 230, 2003, 367-415,
\bibitem[Ra2]{R2}A. Ram, \textit{Seminormal representations of Weyl groups
and Iwahori-Hecke algebras}, Proc. London Math. Soc. (3)  75, no. 1,
99--133, 1997, math.RT/9511223

\bibitem[Rou1]{Rou}R. Rouquier, \textit{Representations of rational
Cherednik algebras}, math.RT/0504600
\bibitem[Rou2]{Rou2}R. Rouquier, \textit{q-Schur algebras and complex
reflection groups, I}, math.RT/0509252

\bibitem[Su1]{Su}T. Suzuki, \textit{Rational and trigonometric
degeneration of the double affine Hecke algebra of type $A$}, IMRN 2005:37
(2005) 2249-2262, math.RT/0502534
\bibitem[Su2]{Su2}T. Suzuki, \textit{Classification of simple modules over
degenerate double affine Hecke algebras of type A},  Int. Math. Res. Not.
2003,  no. 43, 2313--2339,  math.QA/0304474
\bibitem[Su3]{Su3}T. Suzuki, \textit{Double affine Hecke algebras,
conformal coinvariants and Kostka polynomials}, math.QA/0508274

\bibitem[SV]{SV}T. Suzuki, M. Vazirani, \textit{Tableaux on periodic skew
diagrams and irreducible representations of the degenerate double affine
Hecke algebras of type A}, Int. Math. Res. Not.  2005,  no. 27, 1621--1656,
math.QA/0406617

\bibitem[VV1]{VV1} M. Varagnolo, E. Vasserot, \textit{On the decomposition
matrices of the quantized Schur algebra},  Duke Math. J.  100  (1999),
no. 2, 267--297, math.QA/9803023
\bibitem[VV2]{VV2} M. Varagnolo, E. Vasserot, \textit{From double affine
Hecke algebras to quantized affine Schur algebras},  Int. Math. Res. Not.
2004,  no. 26, 1299--1333, math.RT/0307047

\bibitem[Va]{Va}E.Vasserot, \textit{Induced and Simple Modules of Double
Affine Hecke Alegbras}, Duke math. J. 126, 2005, 251-323
\end{thebibliography}
\end{document}